\documentclass[10pt]{article}
\usepackage{amssymb,latexsym}
\usepackage{epsfig}
\usepackage{eufrak}
\usepackage{amsmath}
\usepackage{mathrsfs}
\usepackage{color}

\newtheorem{theorem}{Theorem}
\newtheorem{lemma}{Lemma}

\newcommand{\be}{\begin{equation}}
\newcommand{\ee}{\end{equation}}
\newcommand{\bee}{\begin{eqnarray*}}
\newcommand{\eee}{\end{eqnarray*}}
\newcommand{\bel}{\begin{eqnarray}}
\newcommand{\eel}{\end{eqnarray}}
\newcommand{\bec}{\begin{cases}}
\newcommand{\eec}{\end{cases}}
\newcommand{\bem}{\begin{bmatrix}}
\newcommand{\eem}{\end{bmatrix}}

\newcommand{\la}{\label}
\newcommand{\li}{\left}
\newcommand{\ri}{\right}

\newcommand{\ovl}{\overline}
\newcommand{\udl}{\underline}

\newcommand{\ep}{\epsilon}
\newcommand{\vep}{\varepsilon}
\newcommand{\lm}{\lambda}

\newcommand{\si}{\sigma}

\newcommand{\de}{\delta}

\newcommand{\vse}{\vartheta}
\newcommand{\se}{\theta}

\newcommand{\al}{\alpha}

\newcommand{\ka}{\kappa}

\newcommand{\Om}{\Omega}

\newcommand{\f}{\frac}

\newcommand{\cd}{\cdots}

\newcommand{\qu}{\quad}
\newcommand{\qqu}{\qquad}

\newcommand{\mscr}{\mathscr}
\newcommand{\mcal}{\mathcal}
\newcommand{\mbf}{\mathbf}
\newcommand{\bb}{\mathbb}

\newcommand{\wh}{\widehat}

\newcommand{\mrm}{\mathrm}
\newcommand{\bs}{\boldsymbol}

\newcommand{\tx}{\text}

\newcommand{\iy}{\infty}

\newcommand{\pa}{\partial}

\newcommand{\bed}{\begin{description}}
\newcommand{\eed}{\end{description}}
\newcommand{\bei}{\begin{itemize}}
\newcommand{\eei}{\end{itemize}}
\newcommand{\ben}{\begin{enumerate}}
\newcommand{\een}{\end{enumerate}}
\newcommand{\bib}{\bibitem}
\newcommand{\beL}{\begin{lemma}}
\newcommand{\eeL}{\end{lemma}}
\newcommand{\beT}{\begin{theorem}}
\newcommand{\eeT}{\end{theorem}}
\newcommand{\sect}{\section}

\newcommand{\bpf}{\begin{pf}}
\newcommand{\epf}{\end{pf}}
\newcommand{\bsk}{\bigskip}

\newcommand{\pfbox}{\hfill\mbox{$\Box$}}

\newenvironment{pf}{\paragraph*{Proof{\rm.}}}{\pfbox\bigskip}

\begin{document}

\title{{\bf Sequential Estimation Methods from Inclusion Principle}
\thanks{The author had been previously working with Louisiana
State University at Baton Rouge, LA 70803, USA, and is now with
Department of Electrical Engineering, Southern University and A\&M College,
Baton Rouge, LA 70813, USA; Email: chenxinjia@gmail.com}}

\author{Xinjia Chen}

\date{August 2012}

\maketitle

\begin{abstract}

In this paper, we propose new sequential estimation methods based on inclusion principle. The main idea is to reformulate the estimation
problems as constructing sequential random intervals  and use confidence sequences to control the associated coverage probabilities. In contrast
to existing asymptotic sequential methods, our estimation procedures rigorously guarantee the pre-specified levels of confidence.
\end{abstract}

% \tableofcontents

\section{Introduction}

An important issue of parameter estimation is the determination of sample sizes.  However, the appropriate sample size usually depends on the
parameters to be estimated from the sampling process.  To overcome this difficulty, an adaptive approach, referred to as sequential estimation
method, has been proposed in sequential analysis, where the sample size is not fixed in advance. Instead, data is evaluated as it is collected
and further sampling is stopped in accordance with a pre-defined stopping rule as significant results are observed.  In the area of sequential
estimation, a wide variety of sampling schemes have been proposed to achieve prescribed levels of accuracy and confidence for the estimation
results. Unfortunately, existing sequential estimation methods are dominantly of asymptotic nature.  That is, the guarantee of the pre-specified
confidence level comes only as the margin of error approaches zero or equivalently the average sample size tends to  infinity.  Since any
practical sampling scheme must employ a finite sample size, the application of asymptotic sequential methods inevitably  introduce unknown
statistical error.  To overcome the limitations of existing asymptotic sequential estimation methods, we shall develop new sampling schemes by
virtue of the {\it inclusion principle} proposed in \cite{Chenestimation, Chen_rule}.

In our paper \cite{Chenestimation, Chen_rule}, we have demonstrated that a wide variety of sequential estimation problems can be cast into the
general framework of constructing a {\it sequential random interval} of a prescribed level of coverage probability. To ensure the requirement of
coverage probability, we propose to use a sequence of confidence intervals, referred to as {\it controlling confidence sequence}, to define a
stopping rule such that the sequential random interval must include the controlling confidence sequence at the termination of the sampling
process. In situations that no other requirement imposed on the sequential random interval except the specification of coverage probability, we
have proposed a more specific version of this principle for constructing sampling schemes as follows: {\it The sampling process is continued
until the controlling confidence sequence is included by the sequential random interval at some stage.}  Such a general method of constructing
sequential estimation procedures is referred to as {\it Inclusion Principle}, which can be justified by the following probabilistic results.

\beT \la{General Inclusion Principle} Let $(\Om, \mscr{F}, \{ \mscr{F}_\ell \}, \Pr )$ be a filtered space.   Let $\bs{\tau}$ be a proper
stopping time with support $I_{\bs{\tau}}$.
  For $\ell \in I_{\bs{\tau}}$, let $\bs{\mcal{A}}_\ell$ and
$\bs{\mcal{B}}_\ell$ be random intervals defined by random variables measurable in $\mscr{F}_\ell$.  Assume that {\small $\{ \bs{\tau} = \ell \}
\subseteq \{ \bs{\mcal{A}}_\ell \subseteq \bs{\mcal{B}}_\ell \}$ } for $\ell \in I_{\bs{\tau}}$. Then, $\Pr \{ \se \in \bs{\mcal{B}}_{\bs{\tau}}
\} \geq \Pr \{ \se \in \bs{\mcal{A}}_\ell \; \tx{for} \; \ell \in I_{\bs{\tau}} \} \geq 1 - \sum_{\ell \in I_{\bs{\tau}} } \Pr \{ \se \notin
\bs{\mcal{A}}_\ell \}$ for any real number $\se$. \eeT

See \cite{Chenestimation} for a  proof. This theorem implies that the coverage probability of the sequential random interval constructed based
on the inclusion principle is bounded from below by the coverage probability of the controlling confidence sequence.

The remainder of the present paper is organized as follows.  In Section 2, we shall apply the inclusion principle to develop analytic stopping
rules for estimating the parameters of binomial, geometric and Poisson distributions. For wider applications, we address the problem of
estimating a binomial proportion in a more general setting of estimating the mean of a bounded random  variable.  To make the stopping rules as
simple as possible, we have made effort to eliminate the need of computing confidence limits. In Section 3, we further consider the problem of
estimating the mean of a bounded random variable by taking into account the information of sample variance.  Section 4 is the conclusion. The
justification of stopping rules and proofs of theorems are given in Appendices.  The main results of this paper have appeared in our conference
paper \cite{Chencontrol}.

Throughout this paper, we shall use the following notations.  Let ``$A \vee B$'' denote the maximum of $A$ and $B$.   Let $\bb{N}$ denote the
set of positive integers.  Let $\bb{R}$ denote the set of real numbers.   Let $\Pr \{ E \}$ denote the probability of event $E$.  The
expectation of a random variable is denoted by $\bb{E} [ . ]$. The other notations will be made clear as we proceed.

\section{Analytic Stopping Rules} \la{analytic}

In this section, we shall propose various analytic stopping rules for estimating mean values of random variables with pre-specified precision
and confidence levels.   More formally, let $X$ be a random variable with mean $\bb{E} [ X ]$.  The general problem is to estimate $\bb{E} [ X
]$ based on i.i.d. samples $X_1, X_2, \cd$ of $X$ by virtue of sequential sampling.   For $n \in \bb{N}$, let $\ovl{X}_n$ denote the sample mean
of $X$, i.e., $\ovl{X}_n = \f{ \sum_{i=1}^n X_i}{n}$.  When the sampling process is terminated with sample number $\mbf{n}$, the sample mean
$\ovl{X}_{\mbf{n}}$ is taken as an estimate for $\bb{E} [ X ]$.  To describe our stopping rules, we need to introduce some bivariate functions
as follows.

Define function $\mscr{M}_{\mrm{B}} (., .)$ such that
\[
\mscr{M}_{\mrm{B}} (z,\se) = \bec z \ln \f{\se}{z} + (1 - z) \ln \f{1 - \se}{1 - z} & \; \tx{for} \; z \in (0,1), \; \se \in (0, 1),\\
\ln(1-\se) & \; \tx{for} \; z = 0, \; \se \in (0, 1),\\
\ln \se & \; \tx{for} \; z = 1, \; \se \in (0, 1),\\
- \iy & \; \tx{for} \; z \in (-\iy, \iy), \; \se \notin (0, 1).  \eec
\]
Define function $\mscr{M}_{\mrm{G}} (., .)$ such that
\[
\mscr{M}_{\mrm{G}} (z,\se) = \bec z \ln \f{z}{\se} + \li ( 1 - z \ri ) \ln \f{1 - z}{1 - \se} & \; \tx{for} \; z \in (1, \iy), \; \se \in (1,
\iy),\\
- \ln \se & \; \tx{for} \; z = 1, \; \se \in (1, \iy),\\
- \iy & \; \tx{for} \; z \in [1, \iy), \; \se \notin (1, \iy). \eec
\]
Define function $\mscr{M}_{\mrm{P}} (., .)$ such that
\[
\mscr{M}_{\mrm{P}} (z, \se) = \bec z - \se + z \ln \li ( \f{\se}{z} \ri ) & \; \tx{for} \; z
> 0,  \;  \se > 0,\\
- \se & \; \tx{for} \; z = 0, \; \se > 0,\\
- \iy & \; \tx{for} \; z \geq 0, \; \se \leq 0. \eec
\]

\subsection{Estimation of Means of Bounded Random Variables} \la{EMBRV}

Let $X$ be a random variable such that $\bb{E} [X] = \mu \in (0, 1)$ and $0 \leq X \leq 1$.  Let $\de \in (0, 1)$. To estimate $\mu$ with a
margin of absolute error $\vep \in (0, \f{1}{2})$ and confidence level $1 - \de$, we consider sampling procedures of $s$ stages.  Let $m_1, m_2,
\cd, m_s$ be an ascending sequence of positive integers such that $m_s \geq \f{ \ln \f{2 s}{\de} }{ 2 \vep^2 }$. Let $\mscr{N}$ be a subset of
positive integers which contains $\{ m_1, \cd, m_s \}$. We propose two stopping rules as follows:

\noindent {\bf Stopping Rule A}:  Continue sampling until $\mscr{M}_{\mrm{B}} \li ( \f{1}{2} - \li | \f{1}{2} - \ovl{X}_n \ri | + \vep - \f{n
\vep }{n \vee m_\ell}, \; \f{1}{2} - \li | \f{1}{2} - \ovl{X}_n \ri | + \vep \ri ) \leq \f{1}{m_\ell} \ln \f{\de}{2 s}$ for some integers $n \in
\mscr{N}$ and $\ell \in \{1, \cd, s \}$.

\noindent {\bf Stopping Rule B}: Continue sampling until $\li [ \li | \ovl{X}_n  - \f{1}{2} \ri | - \vep + \f{n \vep }{3 ( n \vee m_\ell ) } \ri
]^2 \geq \f{1}{4} - \li ( \f{n}{n \vee m_\ell} \ri )^2 \f{ m_\ell \vep^2} { 2 \ln \f{2 s}{\de} }$ for some integers $n \in \mscr{N}$ and $\ell
\in \{1, \cd, s \}$.

In Appendices \ref{De_app_A} and \ref{De_app_B}, we have shown that for both stopping rules A and B, the sample mean $\ovl{X}_{\mbf{n}}$ at the
termination of the sampling process guarantees that $\Pr \{ | \ovl{X}_{\mbf{n}} - \mu | < \vep \} \geq 1 - \de$.  To avoid  unnecessary checking
of the stopping conditions,  we suggest choosing $m_1 \geq \f{ \ln \f{2 s}{\de}  }{  \ln \f{1}{1 - \vep}  }$ for Stopping Rule A and $m_1 \geq
\li ( \f{ 24 \vep - 16 \vep^2 }{ 9 }  \ri )  \f{ \ln \f{2 s}{\de} }{ 2 \vep^2 }$ for Stopping Rule B, respectively.  For purpose of efficiency,
we recommend choosing $m_1, \cd, m_s$ as a geometric sequence, i.e., $\f{ m_{\ell + 1} }{m_\ell}$ is approximately equal for $\ell = 1, \cd, s -
1$.

Next, consider the problem of estimating $\mu$ with  a margin of relative error $\vep \in (0, 1)$ and confidence level $1 - \de$.  Let $\de_1,
\de_2, \cd$ be a sequence of positive numbers such that $\sum_{\ell = 1}^\iy \de_\ell = \de \in (0, 1)$.   Let $m_1,  m_2 ,  \cd$ be an
ascending sequence of positive integers such that $\lim_{\ell \to \iy} \f{\ln (\de_\ell)}{m_\ell} = 0$. Let $\mscr{N}$ be a subset of positive
integers which contains $\{ m_1, m_2, \cd \}$.  We propose a stopping rule as follows:

\noindent{\bf Stopping Rule C}: Continue sampling until $\ovl{X}_n > 0$ and $\mscr{M}_{\mrm{B}} \li ( \f{\ovl{X}_n}{1 + \vep} \li ( 1 + \f{n
\vep}{n \vee m_\ell } \ri ), \; \f{\ovl{X}_n}{1 + \vep} \ri ) \leq \f{1}{m_\ell} \ln \f{\de_\ell}{2}$ for some integers $n \in \mscr{N}$ and
$\ell \in \bb{N}$.

In Appendix \ref{De_app_C},  we have established that for stopping rule C, the sampling process will eventually stop with probability $1$ and
the sample mean $\ovl{X}_{\mbf{n}}$ at the termination of the sampling process guarantees that $\Pr \{ | \ovl{X}_{\mbf{n}} - \mu | < \vep \mu \}
\geq 1 - \de$.

\subsection{Estimation of Means of Geometric Distributions} \la{EMGD}

Let $X$ be a random variable having a geometric distribution with mean $\se \in (1, \iy)$.  Let $\vep, \; \de \in (0, 1)$. To estimate $\se$, we
consider sampling procedures of $s$ stages.   Let $m_1, m_2, \cd, m_s$ be an ascending sequence of positive integers. Let $\mscr{N}$ be a subset
of positive integers which contains $\{ m_1, \cd, m_s \}$.  Under the assumption that $m_s \geq \f{(1 + \vep) \ln \f{2 s}{\de} }{ (1 + \vep) \ln
(1 + \vep) - \vep }$, we propose the following stopping rule:

\noindent{\bf Stopping Rule D}: Continue sampling until $\mscr{M}_{\mrm{G}} \li ( \li ( 1 + \vep - \f{n \vep}{n \vee m_\ell} \ri ) \ovl{X}_n, \;
(1 + \vep) \ovl{X}_n \ri ) \leq \f{1}{m_\ell} \ln \f{\de}{2 s}$ for some $n \in \mscr{N}$ and $\ell \in \{1, \cd, s \}$.

In Appendix \ref{De_app_D},  we have proved that for stopping rule D, the sample mean $\ovl{X}_{\mbf{n}}$ at the termination of the sampling
process guarantees that $\Pr \{ (1 - \vep) \ovl{X}_{\mbf{n}} < \se < (1 + \vep) \ovl{X}_{\mbf{n}}  \} \geq 1 - \de$.  To avoid  unnecessary
checking of the stopping condition,  we suggest choosing $m_1 \geq \f{ \ln \f{2 s}{\de}  }{  \ln (1 + \vep) }$ for Stopping Rule D.  For purpose
of efficiency, we recommend choosing $m_1, \cd, m_s$ as a geometric sequence.  It should be noted that the estimation of a binomial proportion
$p$ with a margin of relative error $\vep$ can be accomplished by such method if $\f{1}{p}$ is identified as $\se$.

\subsection{Estimation of Poisson Parameters}

Let $X$ be a random variable having a Poisson distribution with mean $\lm \in (0, \iy)$.  Let $\vep > 0$ and $0 < \de < 1$. To estimate $\lm$,
we consider sampling procedures of infinitely many stages.    Let $\de_1, \de_2, \cd$ be a sequence of positive numbers such that $\sum_{\ell =
1}^\iy \de_\ell = \de \in (0, 1)$.   Let $m_1,  m_2 ,  \cd$ be an ascending sequence of positive integers such that $\lim_{\ell \to \iy} \f{\ln
(\de_\ell)}{m_\ell} = 0$.  Let $\mscr{N}$ be a subset of positive integers which contains $\{ m_1, m_2, \cd \}$.

To estimate $\lm$ with a margin of absolute error $\vep$ and confidence level $1 - \de$, we propose the following stopping rule:

\noindent{\bf Stopping Rule E}:  Continue sampling until $\mscr{M}_{\mrm{P}} \li ( \ovl{X}_n + \vep - \f{n \vep}{n \vee m_\ell}, \; \ovl{X}_n +
\vep \ri ) \leq \f{1}{m_\ell} \ln \f{\de_\ell}{2}$ for some integers $n \in \mscr{N}$ and $\ell \in \bb{N}$.

In Appendix \ref{De_app_E},  we have established that for stopping rule E, the sampling process will eventually stop with probability $1$ and
the sample mean $\ovl{X}_{\mbf{n}}$ at the termination of the sampling process guarantees that $\Pr \{ | \ovl{X}_{\mbf{n}} - \lm | < \vep \}
\geq 1 - \de$.  To estimate $\lm$ with a margin of relative error $\vep$ and confidence level $1 - \de$, we propose the following stopping rule:

\noindent{\bf Stopping Rule F}:  Continue sampling until $\ovl{X}_n > 0$ and $\mscr{M}_{\mrm{P}} \li ( \f{\ovl{X}_n}{1 + \vep} \li ( 1 + \f{n
\vep}{n \vee m_\ell} \ri ), \; \f{\ovl{X}_n}{1 + \vep} \ri ) \leq \f{1}{m_\ell} \ln \f{\de_\ell}{2}$ for some integers $n \in \mscr{N}$ and
$\ell \in \bb{N}$.

In Appendix \ref{De_app_F},  we have established that for stopping rule F, the sampling process will eventually stop with probability $1$ and
the sample mean $\ovl{X}_{\mbf{n}}$ at the termination of the sampling process guarantees that $\Pr \{ | \ovl{X}_{\mbf{n}} - \lm | < \vep \lm \}
\geq 1 - \de$.

\section{Estimation of Means and Variances of Bounded Variables}

In Section \ref{EMBRV},  we have proposed sequential methods for estimating the mean of a bounded random variable.  However, the information of
sample variance is not used in these methods.  In this section, we shall exploit the information of sample variance for purpose of improving the
efficiency of estimation.  To apply the inclusion principle to construct an estimation procedure for estimating the mean of a bounded random
variable, we need to have a confidence sequence for the mean.  The construction of the required confidence sequence can be accomplished by
applying Bonferroni's inequality to a sequence of fixed-sample-size confidence intervals.   Therefore, in the sequel, we shall first study the
construction of fixed-sample-size confidence intervals for the mean and variance of a bounded random variable.

Since any bounded random variable can be expressed as a linear function of a random variable bounded in $[0, 1]$,  it will loss no generality to
consider a random variable $X$  bounded in interval $[0, 1]$, which has mean $\mu \in (0, 1)$ and variance $\si^2 > 0$.  Let $X_1, \cd, X_n$ be
i.i.d. samples of $X$.  Define $\ovl{X} = \f{\sum_{i=1}^n X_i}{n}$ and $\ovl{V} = \f{ \sum_{i = 1}^n (X_i - \ovl{X})^2 }{n}$.  In many
situations, it is desirable to construct confidence intervals for $\mu$ and $\si^2$ based on $\ovl{X}$ and $\ovl{V}$.   For this purpose, we
need to make use of Hoeffding's inequalities.  Specifically, define {\small
\[
\varphi (z, \nu, \se) = \li ( 1 - \f{z \nu }{\nu^2 + \se} \ri ) \ln \f{\se + \nu (\nu - z) } {\se} + \f{z \nu }{\nu^2 + \se} \ln \f{z}{\nu}
\] } for $0 < z < \nu < 1$ and $0 < \se < 1$.  Define $\psi(z, \nu, \se)  = \varphi (1 - z, 1 - \nu, \se)$ for $0 < \nu < z < 1$ and $0 < \se < 1$.
Define $\phi (z, \se) = (1 - z) \ln \f{1 - z}{1 - \se} + z \ln \f{z}{\se}$ for $0 < z < 1$ and $0 < \se < 1$.   Hoeffding's inequalities assert
that \bee &  & \Pr \{ \ovl{X} \geq z \} \leq \exp ( - n \psi(z, \mu, \si^2) ) \leq \exp ( - n \phi(z, \mu)
) \qqu \tx{for $0 < \mu < z$},\\
&   & \Pr \{ \ovl{X} \leq z \} \leq \exp \li ( - n \varphi (z, \mu, \si^2)  \ri ) \leq \exp ( - n \phi(z, \mu) ) \qqu \tx{for $z < \mu < 1$}.
\eee

We have the following results.

\beT \la{seegood}

\bel &  & \f{\pa \psi (z, \mu, \se)}{\pa \mu} \leq 0 \qqu \tx{for $0 < \mu < z$}, \la{good1}\\
&  & \f{\pa \psi (z, \mu, \se)}{\pa \se} \leq 0 \qqu \tx{for $0 < \se < 1$}, \la{good2}\\
&  & \f{\pa \varphi (z, \mu, \se)}{\pa \mu}  \geq 0 \qqu \tx{for $0 < z < \mu$}, \la{good3}\\
&  & \f{\pa \varphi (z, \mu, \se)}{\pa \se} \leq 0 \qqu \tx{for $0 < \se < 1$}. \la{good4}
 \eel \eeT

See Appendix \ref{seegood_app} for a proof.

\subsection{Confidence Interval for Mean Value} \la{CImean}

For simplicity of notations, define $W_\nu = \ovl{V} + (\ovl{X} - \nu)^2$ for $0 \leq \nu \leq 1$.  For constructing a confidence interval for
the mean, we have the following method.

\beT \la{meaninterval} Let $\de \in (0, 1)$.  Define {\small \bee &  & L = \bec \sup \li \{ \nu \in (0, \ovl{X}): \; \max \li [  \psi ( \ovl{X},
\nu, \vse), \; \phi (W_\nu, \vse) \; \bb{I}_{\{\vse > W_\nu\}}  \ri ]  > \f{\ln \f{3}{\de}}{n} \; \tx{for all} \; \vse \in \li (0, \nu (1
- \nu) \ri ] \ri \} & \; \tx{if} \; \ovl{X} > 0,\\
0  & \; \tx{if} \; \ovl{X} = 0
\eec  \\
&   & U = \bec \inf \li \{ \nu \in (\ovl{X}, 1): \; \max \li [  \varphi ( \ovl{X}, \nu, \vse), \; \phi (W_\nu, \vse ) \; \bb{I}_{\{\vse >
W_\nu\}} \ri ]
> \f{\ln \f{3}{\de}}{n} \; \tx{for all} \; \vse \in \li (0, \nu (1 - \nu) \ri ]  \ri \} & \; \tx{if} \; \ovl{X} < 1,\\
1 & \; \tx{if} \; \ovl{X} = 1, \eec   \eee} where $\bb{I}_{\{\vse > W_\nu\}}$ is the indicator function which takes value $1$ if $\vse > W_\nu$
and otherwise tales value $0$.
 Then, $\Pr \{ L \leq \mu \leq U \} \geq 1 - \de$.   \eeT

See Appendix  \ref{meaninterval_app} for a proof.  The computation of the confidence limits is addressed in the sequel.

\subsubsection{Adaptive Scanning Algorithms}

To compute the lower confidence limit $L$, we first need to establish a method to check, for a given interval $[a, b] \subseteq [0, \ovl{X}]$,
whether the following statement is true: \be \la{stateA} \tx{ For every $\nu \in [a, b]$, $\max \li [ \psi ( \ovl{X}, \nu, \vse), \; \phi
(W_\nu, \vse) \; \bb{I}_{\{\vse > W_\nu\}}  \ri ]  > \f{\ln \f{3}{\de}}{n} \; \tx{for all} \; \vse \in \li (0, \nu (1 - \nu) \ri ]$ }. \ee To
check the truth of (\ref{stateA}) without exhaustive computation, our approach is to find a sufficient condition for (\ref{stateA}) so that the
conservativeness of the sufficient condition diminishes as the width of the interval $[a, b]$ decreases.   For simplicity of notations, let $c =
\max \{ a (1 - a), b (1 - b) \}$.  As a consequence of (\ref{good1}), we have $\psi ( \ovl{X}, \nu, \vse) \geq \psi ( \ovl{X}, b, \vse)$ for all
$\nu \in [a, b]$.  Since $W_b \leq W_\nu \leq W_a$ for $\nu \in [a, b]$ and $\phi (z, \vse)$ is non-increasing with respect to $z \in (0,
\vse)$,  we have $\phi (W_\nu, \vse) \bb{I}_{\{\vse > W_\nu\}} \geq \phi (W_a, \vse) \bb{I}_{\{\vse > W_a \}}$ for $\nu \in [a, b]$. Hence, a
sufficient condition for (\ref{stateA}) is as follows: \be \la{stateB} \max \li [ \psi ( \ovl{X}, b, \vse), \; \phi (W_a, \vse) \;
\bb{I}_{\{\vse > W_a \}}  \ri ]  > \f{\ln \f{3}{\de}}{n} \; \tx{for all} \; \vse \in \li (0, c \ri ]. \ee The truth of statement (\ref{stateB})
can be checked by virtue of the following facts:

\bei

\item In the case of  $W_a \geq c$, it follows from (\ref{good2}) that statement (\ref{stateB}) is true if and only if  $\psi ( \ovl{X}, b, c) > \f{\ln \f{3}{\de}}{n}$.

\item In the case of  $W_a < c$, it follows from (\ref{good2}) that statement (\ref{stateB})  is true if and only if
\bel &  & \psi ( \ovl{X}, b, W_a) > \f{\ln
  \f{3}{\de}}{n}, \nonumber \\
  &  & \max \li [ \psi ( \ovl{X}, b, \vse), \; \phi (W_a, \vse ) \ri ] >   \f{\ln
  \f{3}{\de}}{n} \qu \tx{for all} \; \vse \in \li ( W_a, c \ri ].  \la{stateB2} \eel

\eei The truth of statement (\ref{stateB2}) can be checked by making use of the following observations:

\bei

\item In the case of  $\phi (W_a, c ) \leq \f{\ln \f{3}{\de}}{n}$, we have $\phi (W_a, \vse ) \leq \f{\ln \f{3}{\de}}{n}$ for all $\vse \in \li ( W_a, c \ri
]$, since $ \phi (W_a, \vse)$ is non-decreasing with respect to $\vse \in (W_a, c)$.  It follows from (\ref{good2}) that statement
(\ref{stateB2}) is true if and only if $\psi ( \ovl{X}, b, c) > \f{\ln \f{3}{\de}}{n}$.

\item In the case of $\phi (W_a, c ) > \f{\ln \f{3}{\de}}{n}$, there exists a $\se^* \in (W_a, c)$ such that $\phi (W_a, \se^* ) = \f{\ln
\f{3}{\de}}{n}$, since $\phi (W_a, \se)$ is non-decreasing with respect to $\se \in (W_a, c)$.  Thus, $\phi (W_a, \vse) \leq \f{\ln
\f{3}{\de}}{n}$ for all $\vse \in \li ( W_a, \se^* \ri ]$. It follows from (\ref{good2}) that statement (\ref{stateB2}) is true if and only if
$\psi ( \ovl{X}, b, \se^* ) > \f{\ln \f{3}{\de}}{n}$. In practice, $\se^*$ can be replaced by a lower bound $\udl{\se}$ which is extremely tight
(for example, $0 < \se^* - \udl{\se} < 10^{-10}$).  Such a lower bound $\udl{\se}$ can be obtained by a bisection search method.

\eei

\bsk

Therefore, through the above discussion, we have developed a rigorous  method for checking the truth of (\ref{stateB}).  Based on this critical
subroutine,  we propose an efficient method for computing the lower confidence limit $L$ for $\ovl{X} > 0$  as follows.

\bsk

\begin{tabular} {|l |}
\hline
$ \nabla  \; \tx{Choose initial step size $d > \eta$, where $\eta$ is an extremely small number(e.g., $10^{-15}$) }$.\\
$ \nabla \; \tx{Let $F \leftarrow 0$ and $a \leftarrow 0$}$.\\
$ \nabla \; \tx{While $F = 0$, do the following}$:\\
$ \indent \indent \; \diamond  \; \tx{Let $\tx{st} \leftarrow 0$ and $\ell \leftarrow 2$}$;\\
$ \indent \indent \; \diamond \; \tx{While $\tx{st} = 0$, do the following}$:\\
$\indent \indent \indent \indent \; \star \; \tx{Let $\ell \leftarrow \ell - 1$ and $d \leftarrow d  2^\ell$}$.\\
$ \indent \indent \indent \indent \; \star \; \tx{If $a + d < \ovl{X}$, then let $b \leftarrow a + d$.
If (\ref{stateB}) holds, then let $\tx{st} \leftarrow 1$ and $a \leftarrow b$}$.\\
$ \indent \indent \indent \indent \; \star \; \tx{If $d < \eta$, then let $\tx{st} \leftarrow 1$ and $F \leftarrow 1$}$.\\
$ \nabla \; \tx{Return $a$ as the lower confidence limit $L$ for $\ovl{X} > 0$}$.
\\ \hline
\end{tabular}

\bsk

We call this algorithm as {\it Adaptive Scanning Algorithm}, since
it adaptively scans the interval $[0, \ovl{X}]$ to check the truth
of (\ref{stateB}).

To compute the upper confidence limit $U$, we first need to
establish a method to check, for a given interval $[a, b] \subseteq
[\ovl{X}, 1]$, whether the following statement is true: \be
\la{stateAU} \tx{ For every $\nu \in [a, b]$, $\max \li [ \varphi (
\ovl{X}, \nu, \vse), \; \phi (W_\nu, \vse) \; \bb{I}_{\{\vse >
W_\nu\}}  \ri ]  > \f{\ln \f{3}{\de}}{n} \; \tx{for all} \; \vse \in
\li (0, \nu (1 - \nu) \ri ]$ }. \ee To check the truth of
(\ref{stateAU}) without exhaustive computation, our approach is to
find a sufficient condition for (\ref{stateAU}) so that the
conservativeness of the sufficient condition diminishes as the width
of the interval $[a, b]$ decreases.   For simplicity of notations,
let $c = \max \{ a (1 - a), b (1 - b) \}$ as before.  As a
consequence of (\ref{good3}), we have $\varphi ( \ovl{X}, \nu, \vse)
\geq \varphi ( \ovl{X}, a, \vse)$ for all $\nu \in [a, b]$.  Since
$W_a \leq W_\nu \leq W_b$ for $\nu \in [a, b]$ and $\phi (z, \vse)$
is non-increasing with respect to $z \in (0, \vse)$,  we have $\phi
(W_\nu, \vse ) \bb{I}_{\{\vse > W_\nu\}} \geq \phi (W_b, \vse)
\bb{I}_{\{\vse > W_b \}}$ for $\nu \in [a, b]$. Hence, a sufficient
condition for (\ref{stateAU}) is as follows: \be \la{stateBU} \max
\li [ \varphi ( \ovl{X}, a, \vse), \; \phi (W_b, \vse ) \;
\bb{I}_{\{\vse
> W_b \}}  \ri ]  > \f{\ln \f{3}{\de}}{n} \; \tx{for all} \; \vse \in \li (0, c \ri ]. \ee The truth of statement (\ref{stateBU}) can be checked
by virtue of the following facts:

\bei

\item In the case of  $W_b \geq c$, it follows from (\ref{good4}) that statement (\ref{stateBU}) is true if and only if  $\varphi ( \ovl{X}, a, c) > \f{\ln \f{3}{\de}}{n}$.

\item In the case of  $W_b < c$,  it follows from (\ref{good4}) that statement (\ref{stateBU})  is true if and only if
\bel &  & \varphi ( \ovl{X}, a, W_b) > \f{\ln
  \f{3}{\de}}{n}, \nonumber \\
  &  & \max \li [ \varphi ( \ovl{X}, a, \vse), \; \phi (W_b, \vse ) \ri ] >   \f{\ln
  \f{3}{\de}}{n} \qu \tx{for all} \; \vse \in \li ( W_b, c \ri ].  \la{stateBU2} \eel

\eei The truth of statement (\ref{stateBU2}) can be checked by making use of the following observations:

\bei

\item In the case of  $\phi (W_b, c ) \leq \f{\ln \f{3}{\de}}{n}$,
we have $\phi (W_b, \vse ) \leq \f{\ln \f{3}{\de}}{n}$ for all $\vse
\in \li ( W_b, c \ri ]$, since $ \phi (W_b, \vse)$ is non-decreasing
with respect to $\vse \in (W_b, c)$.  It follows from (\ref{good4})
that statement (\ref{stateBU2}) is true if and only if $\varphi (
\ovl{X}, a, c) > \f{\ln \f{3}{\de}}{n}$.

\item In the case of $\phi (W_b, c ) > \f{\ln \f{3}{\de}}{n}$, there exists a $\se^\star \in (W_b, c)$ such that $\phi (W_b, \se^\star) = \f{\ln
\f{3}{\de}}{n}$, since $\phi (W_b, \se)$ is non-decreasing with respect to $\se \in (W_b, c)$.    Thus, $\phi (W_b, \vse) \leq \f{\ln
\f{3}{\de}}{n}$ for all $\vse \in \li ( W_b, \se^\star \ri ]$.  It follows from (\ref{good4}) that statement (\ref{stateBU2}) is true if and
only if $\varphi ( \ovl{X}, a, \se^\star) \leq \f{\ln \f{3}{\de}}{n}$.   In practice, $\se^\star$ can be replaced by a lower bound $\udl{\se}$
which is extremely tight (for example,  $0 < \se^\star - \udl{\se} < 10^{-10}$).  Such a lower bound $\udl{\se}$ can be obtained by a bisection
search method.

\eei

Therefore, through the above discussion, we have developed a rigorous  method for checking the truth of (\ref{stateBU}).  Based on this critical
subroutine,  we propose an efficient method for computing the upper confidence limit $U$ for $\ovl{X} < 1$ as follows. \bsk

\begin{tabular} {|l |}
\hline
$ \nabla  \; \tx{Choose initial step size $d > \eta$, where $\eta$ is an extremely small number(e.g., $10^{-15}$) }$.\\
$ \nabla \; \tx{Let $F \leftarrow 0$ and $b \leftarrow 1$}$.\\
$ \nabla \; \tx{While $F = 0$, do the following}$:\\
$ \indent \indent \; \diamond  \; \tx{Let $\tx{st} \leftarrow 0$ and $\ell \leftarrow 2$}$;\\
$ \indent \indent \; \diamond \; \tx{While $\tx{st} = 0$, do the following}$:\\
$\indent \indent \indent \indent \; \star \; \tx{Let $\ell \leftarrow \ell - 1$ and $d \leftarrow d  2^\ell$}$.\\
$ \indent \indent \indent \indent \; \star \; \tx{If $b - d > \ovl{X}$, then let $a \leftarrow b - d$.
If (\ref{stateBU}) holds, then let $\tx{st} \leftarrow 1$ and $b \leftarrow a$}$.\\
$ \indent \indent \indent \indent \; \star \; \tx{If $d < \eta$, then let $\tx{st} \leftarrow 1$ and $F \leftarrow 1$}$.\\
$ \nabla \; \tx{Return $b$ as the upper confidence limit $U$ for $\ovl{X} < 1$}$.
\\ \hline
\end{tabular}

\bsk

We call this algorithm as {\it Adaptive Scanning Algorithm}, since
it adaptively scans the interval $[\ovl{X}, 1]$ to check the truth
of (\ref{stateBU}).

\subsection{Sequential Estimation of Mean}

In the preceding discussion, we have developed rigorous methods for constructing fixed-sample-size confidence intervals for the mean $\mu$ of
the random variable $X$ bounded in $[0, 1]$.  Now, we are ready to construct a multistage sampling scheme which produces an estimator
$\wh{\bs{\mu}}$ for $\mu$ such that $\Pr \{  | \wh{\bs{\mu}} - \mu | < \vep  \} \geq 1 - \de$,  where $\vep, \; \de \in (0, 1)$.   For this
purpose, we consider a sampling procedure of $s$ stages, with sample sizes $n_1 < n_2 < \cd < n_s$ chosen such that
\[
\f{ \ln \f{2 s}{\de}  }{  \ln \f{1}{1 - \vep}  } \leq  n_1 \leq \f{ \ln \f{2 s}{\de}  }{ 2 \vep^2 } \leq n_s.
\]
At each stage with index $\ell \in \{1, \cd, s \}$, we use the method described in Section \ref{CImean} to construct a confidence interval
$(L_\ell, U_\ell)$ for $\mu$ in terms of $\ovl{X}_{n_\ell} = \f{ \sum_{i=1}^{n_\ell} X_i }{n_\ell}$ and $\ovl{V}_{n_\ell}  = \f{ \sum_{i = 1}^n
(X_i - \ovl{X}_{n_\ell})^2 }{n_\ell}$ such that $\Pr \{ L_\ell < \mu < U_\ell \} \geq 1 - \f{\de}{2 s}$.  Then, from Bonferroni's inequality, we
have a confidence sequence $\{ (L_\ell, U_\ell), \; 1 \leq \ell \leq s \}$ such that $\Pr \{ L_\ell < \mu < U_\ell, \; \ell = 1, \cd, s \} \geq
1 - \de$.  By the inclusion principle, a stopping rule can be defined as follows:

Continue sampling until there exists an index $\ell \in \{ 1, \cd, s \}$ such that $\ovl{X}_{n_\ell} - \vep \leq L_\ell \leq U_\ell \leq
\ovl{X}_{n_\ell} + \vep$.  At the termination of the sampling process, take $\ovl{X}_{n_\ell}$ with the corresponding index $\ell$ as the
estimator $\wh{\bs{\mu}}$ for $\mu$.

According to Theorem \ref{General Inclusion Principle}, the estimator $\wh{\bs{\mu}}$ for $\mu$ resulted from the above procedure ensures that
$\Pr \{  | \wh{\bs{\mu}} - \mu | < \vep  \} \geq 1 - \de$.

\subsection{Confidence Region for Mean and Variance}

In many situations, it might be interested to infer both the mean $\mu$ and variance $\si^2$ of $X$.  For constructing confidence region for the
mean $\mu$ and variance $\si^2$, we propose the following method.

\beT \la{Region} Let $\de \in (0, 1)$.  Define \bee &  & \mscr{A} = \li \{ (\nu, \vse):  \ovl{X} \leq \nu < 1, \; 0 < \vse \leq \nu (1 - \nu),
\; \varphi ( \ovl{X}, \nu, \vse) <
\f{1}{n} \ln \f{4}{\de}, \; \phi (W_\nu, \vse) < \f{1}{n} \ln \f{4}{\de}  \ri \},\\
&  & \mscr{B} = \li \{ (\nu, \vse):  \ovl{X} > \nu > 0, \; 0 < \vse \leq \nu (1 - \nu), \;  \psi ( \ovl{X}, \nu, \vse) < \f{1}{n} \ln
\f{4}{\de}, \; \phi (W_\nu, \vse) < \f{1}{n} \ln \f{4}{\de} \ri  \} \eee and $\mscr{D} ( \ovl{X}, \ovl{V}) = \mscr{A} \cup \mscr{B}$. Then, $\Pr
\{ (\mu, \si^2) \in \mscr{D} ( \ovl{X}, \ovl{V} ) \} \geq 1 - \de$.   \eeT

See Appendix  \ref{Region_app} for a proof.  The boundary of $\mscr{A}$ is a subset of $C_1 \cup C_2 \cup C_3$,
where \bee &   & C_1 = \li \{ (\nu, \vse):  \ovl{X} \leq \nu < 1, \; \vse = \nu (1 - \nu) \ri \}, \\
&   & C_2 = \li \{ (\nu, \vse): \ovl{X} \leq \nu < 1, \; 0 < \vse \leq \f{1}{4}, \; \varphi ( \ovl{X}, \nu, \vse) =
\f{1}{n} \ln \f{4}{\de} \ri \},\\
&  & C_3 = \li \{ (\nu, \vse):  \ovl{X} \leq \nu < 1, \; W_\nu < \vse \leq \f{1}{4}, \; \phi (W_\nu, \vse) = \f{1}{n} \ln \f{4}{\de}  \ri \}\\
&  & \qqu \qu \bigcup \li \{ (\nu, \vse):  \ovl{X} \leq \nu < 1, \; 0 < \vse < W_\nu, \; \phi (W_\nu, \vse) = \f{1}{n} \ln \f{4}{\de} \ri \}.
\eee As a consequence of (\ref{good3}), $\varphi ( \ovl{X}, \nu, \vse)$ is non-decreasing with respect to $\nu$.  Hence, the points in $C_2$ can
be obtained by solving equation $\varphi ( \ovl{X}, \nu, \vse) = \f{1}{n} \ln \f{4}{\de}$ for $\nu \in [\ovl{X}, 1)$ with a bisection search
method.  Note that $\phi (W_\nu, \vse)$ is non-increasing with respect to $\vse \in (0, W_\nu)$ and is non-decreasing with respect to $\vse \in
(W_\nu, \f{1}{4})$.  It follows that the points in $C_3$ can be obtained by solving equation $\phi (W_\nu, \vse) =  \f{1}{n} \ln \f{4}{\de}$ for
$\vse$ with a bisection search method.

On the other side, the boundary of $\mscr{B}$ is a subset of $D_1 \cup D_2 \cup D_3$,
where \bee &   & D_1 = \li \{ (\nu, \vse):  \ovl{X} > \nu > 0, \; \vse = \nu (1 - \nu) \ri \}, \\
&   & D_2 = \li \{ (\nu, \vse): \ovl{X} > \nu > 0, \; 0 < \vse \leq \f{1}{4}, \; \psi ( \ovl{X}, \nu, \vse) =
\f{1}{n} \ln \f{4}{\de} \ri \},\\
&  & D_3 = \li \{ (\nu, \vse):  \ovl{X} > \nu > 0, \; W_\nu < \vse \leq \f{1}{4}, \; \phi (W_\nu, \vse) = \f{1}{n} \ln \f{4}{\de}  \ri \}\\
&  & \qqu \qu \bigcup \li \{ (\nu, \vse):  \ovl{X} > \nu > 0, \; 0 < \vse < W_\nu, \; \phi (W_\nu, \vse) = \f{1}{n} \ln \f{4}{\de} \ri \}. \eee
As a consequence of (\ref{good1}), $\psi ( \ovl{X}, \nu, \vse)$ is non-increasing with respect to $\nu$.  Hence, the points in $D_2$ can be
obtained by solving equation $\psi ( \ovl{X}, \nu, \vse) = \f{1}{n} \ln \f{4}{\de}$ for $\nu \in [\ovl{X}, 1)$ with a bisection search method.
Note that $\phi (W_\nu, \vse)$ is non-increasing with respect to $\vse \in (0, W_\nu)$ and is non-decreasing with respect to $\vse \in (W_\nu,
\f{1}{4})$.  It follows that the points in $D_3$ can be obtained by solving equation $\phi (W_\nu, \vse) =  \f{1}{n} \ln \f{4}{\de}$ for $\vse$
with a bisection search method.

Finally, we would like to point out that one can apply the same technique to develop confidence intervals and sequential estimation procedures
for the mean and variance  based on bounding the tail probabilities $\Pr \li \{ \ovl{X} \geq z \ri \}$ and $\Pr \li \{ \ovl{X} \leq z \ri \}$ by
Bennet's inequalities \cite{Bennett} or Bernstein's inequalities \cite{Bernstein}.

\sect{Conclusion}

In this paper, we have applied inclusion principle to develop extremely simple analytic sequential methods for estimating the means of binomial,
geometric, Poisson and bounded random variables.  Moreover,  we have developed sequential methods for estimating the mean of bounded random
variables, which makes use of the information of sample variance.  Our sequential estimation methods guarantee the prescribed levels of accuracy
and confidence.

\appendix

\sect{Derivation of Stopping Rules}

For simplicity of notations, define $\mscr{S} = \{1, \cd, s \}$.

\subsection{Derivation of Stopping Rule A} \la{De_app_A}

We need some preliminary results.  As applications of Corollary 5 of
\cite{Chenmax}, we have Lemmas \ref{lem1338} and \ref{lem1338b}.

\beL \la{lem1338}

Let $\mu \in (0, 1)$.  Let $m \in \mscr{N}$ and $\vep \in (0, 1 - \mu)$.  Then,
\[
\Pr \li \{  \ovl{X}_n < \mu + \f{ (m \vee n) \vep }{n} \; \tx{for all} \; n \in \mscr{N} \ri \} \geq 1 - \exp \li ( m \mscr{M}_{\mrm{B}} ( \mu +
\vep, \mu ) \ri ).
\]
\eeL

\beL \la{lem1338b}

Let $\mu \in (0, 1)$.  Let $m \in \mscr{N}$ and $\vep \in (0, \mu)$.  Then,
\[
\Pr \li \{  \ovl{X}_n > \mu - \f{(m \vee n) \vep }{n} \; \tx{for all} \; n \in \mscr{N} \ri \} \geq 1 - \exp \li ( m \mscr{M}_{\mrm{B}} ( \mu -
\vep, \mu ) \ri ).
\]
\eeL

\beL   \la{lem2338} Let $y, \; r \in (0, 1]$.  Then, $\mscr{M}_{\mrm{B}} \li ( \mu + r \li (  y - \mu \ri ), \mu \ri )$ is non-decreasing  with
respect to $\mu \in (0, y)$.

\eeL

\bpf  From the definition of the function $\mscr{M}_{\mrm{B}}$, we have that $\mscr{M}_{\mrm{B}} ( z, \mu) = z \ln \f{\mu}{z} + (1 - z) \ln \f{1
- \mu}{1 - z}$ for $z \in (0, 1)$ and $\mu \in (0, 1)$.   It can be checked that $\f{ \pa  \mscr{M}_{\mrm{B}} ( z, \mu) }{ \pa z } =  \ln \f{\mu
(1 - z)}{z (1 - \mu)}$ and $\f{ \pa  \mscr{M}_{\mrm{B}} ( z, \mu) }{ \pa \mu } = \f{z - \mu}{\mu (1 - \mu)}$ for $z \in (0, 1)$ and $\mu \in (0,
1)$.  Now let $z =  \mu + r \li ( y - \mu \ri )$.   Since $\mu \in (0, y)$, it follows that $z \in (0, 1)$.   Hence, \bee  \f{ \pa
\mscr{M}_{\mrm{B}} \li ( \mu + r \li ( y - \mu \ri ), \mu \ri )}{\pa \mu} & = & ( 1 - r)  \ln \f{\mu (1
- z)}{z (1 - \mu)} + \f{z - \mu}{\mu (1 - \mu)} \\
& \geq  &  \f{( 1 - r) (\mu - z) }{\mu (1 - z)} -  \f{\mu - z}{\mu (1 - \mu)}  =  \f{ r (\mu - z) (y - 1) }{ \mu (1 - \mu) ( 1 - z )} \geq 0.
\eee This completes the proof of the lemma.

\epf

\beL   \la{lem2338b} Let $y \in [0, 1)$ and $r \in (0, 1]$.  Then, $\mscr{M}_{\mrm{B}} \li ( \mu - r \li (  \mu - y \ri ), \mu \ri )$ is
non-increasing  with respect to $\mu \in (y, 1)$.

\eeL

\bpf

For simplicity of notations, let $z = \mu - r \li (  \mu - y \ri )$.   Note that \bee \f{ \pa \mscr{M}_{\mrm{B}} \li ( \mu - r \li (  \mu - y
\ri ), \mu \ri )}{\pa \mu} & = & ( 1 - r) \ln \f{\mu (1 -
z)}{z (1 - \mu)} + \f{z - \mu}{\mu (1 - \mu)}\\
& \leq &  \f{( 1 - r) (\mu - z) }{z (1 - \mu)}  -  \f{\mu - z}{ \mu (1 - \mu)}  =  \f{ r (z - \mu) y } { \mu (1 - \mu) z } \leq 0. \eee  This
proves the lemma.

\epf

\beL \la{lem3338}  For $n \in \mscr{N}$ and $\ell \in \mscr{S}$, define \[ L_n^\ell = \bec \inf \li \{ \nu \in (0, \ovl{X}_n):
\mscr{M}_{\mrm{B}} \li ( \nu + \f{n}{n \vee m_\ell} \li (
\ovl{X}_n - \nu \ri ), \nu \ri ) > \f{1}{m_\ell} \ln \f{\de}{2 s} \ri \}  & \; \tx{for} \; \ovl{X}_n > 0,\\
0 & \; \tx{for} \; \ovl{X}_n = 0.  \eec \]  Then, $\Pr \{ L_n^\ell < \mu \; \tx{for all} \; n \in \mscr{N} \} \geq 1 - \f{\de}{2 s}$ for $\ell
\in \mscr{S}$.  \eeL

\bpf

First, we need to show that $L_n^\ell$ is well-defined.  Since $L_n^\ell = 0$ for $\ovl{X}_n = 0$, $L_n^\ell$ is well-defined provided that
$L_n^\ell$ exists for $0 < \ovl{X}_n \leq 1$.  Note that  $\lim_{\nu \uparrow y} \mscr{M}_{\mrm{B}} ( \nu + \f{n}{n \vee m_\ell} ( y - \nu ),
\nu ) = 0 > \f{1}{m_\ell} \ln \f{\de}{2 s}$ for $y \in (0, 1]$.  This fact together with Lemma \ref{lem2338} imply the existence of $L_n^\ell$
for $0 < \ovl{X}_n \leq 1$. So, $L_n^\ell$ is well-defined.  From the definition of $L_n^\ell$, it can be seen that \bee &  & \{ \mu \leq
L_n^\ell, \; \ovl{X}_n = 0 \} = \li \{ \mu \leq \ovl{X}_n, \;
\mscr{M}_{\mrm{B}} \li ( \mu + \f{n}{n \vee m_\ell} ( \ovl{X}_n - \mu ), \mu \ri ) \leq \f{1}{m_\ell} \ln \f{\de}{2 s}, \; \ovl{X}_n = 0 \ri \} = \emptyset,\\
&  & \{ \mu \leq L_n^\ell, \; 0 < \ovl{X}_n \leq 1 \} \subseteq  \li \{ \mu \leq \ovl{X}_n, \; \mscr{M}_{\mrm{B}} \li ( \mu + \f{n}{n \vee
m_\ell} ( \ovl{X}_n - \mu ), \mu \ri ) \leq \f{1}{m_\ell} \ln \f{\de}{2 s}, \; 0 < \ovl{X}_n \leq 1 \ri \}. \eee This implies that $\{ \mu \leq
L_n^\ell \} \subseteq  \{ \mu \leq \ovl{X}_n, \; \mscr{M}_{\mrm{B}} ( \mu + \f{n}{n \vee m_\ell} ( \ovl{X}_n - \mu ), \mu ) \leq \f{1}{m_\ell}
\ln \f{\de}{2 s} \}$.

Next, consider $\Pr \{ L_n^\ell < \mu \; \tx{for all} \; n \in \mscr{N} \}$ for two cases as follows.

Case A: $\mu^{m_\ell} \leq \f{\de}{2 s}$.

Case B: $\mu^{m_\ell} > \f{\de}{2 s}$.

In Case A, there must exist an $\vep^* \in (0, 1 - \mu]$ such that
$\mscr{M}_{\mrm{B}} \li ( \mu + \vep^*, \mu \ri ) = \f{1}{m_\ell}
\ln \f{\de}{2 s}$. Note that $\mscr{M}_{\mrm{B}} ( \mu + \ep, \mu )$
is decreasing with respect to $\ep \in (0, 1 - \mu)$. Therefore,
from the definitions of $L_n^\ell$ and $\vep^*$, we have that $\{
\mu \leq L_n^\ell \} \subseteq  \{ \mu \leq \ovl{X}_n, \;
\mscr{M}_{\mrm{B}} ( \mu + \f{n}{n \vee m_\ell} ( \ovl{X}_n - \mu ),
\mu ) \leq \f{1}{m_\ell} \ln \f{\de}{2 s} \} \subseteq \{ \mu \leq
\ovl{X}_n, \; \f{n}{n \vee m_\ell} ( \ovl{X}_n - \mu ) \geq \vep^*
\} \subseteq \{ \ovl{X}_n \geq \mu + \f{(n \vee m_\ell) \vep^* }{n}
\}$. This implies that $\{ L_n^\ell < \mu \} \supseteq \{ \ovl{X}_n
< \mu + \f{(n \vee m_\ell) \vep^* }{n} \}$ for all $n \in \mscr{N}$.
Hence, $\{ L_n^\ell < \mu \; \tx{for all} \; n \in \mscr{N} \}
\supseteq \{ \ovl{X}_n < \mu + \f{(n \vee m_\ell) \vep^* }{n} \;
\tx{for all} \; n \in \mscr{N} \}$. It follows from Lemma
\ref{lem1338} that $\Pr \{ L_n^\ell < \mu \; \tx{for all} \; n \in
\mscr{N} \} \geq \Pr \{ \ovl{X}_n < \mu + \f{(n \vee m_\ell) \vep^*
}{n} \; \tx{for all} \; n \in \mscr{N} \} \geq 1 - \exp \li ( m_\ell
\mscr{M}_{\mrm{B}} ( \mu + \vep^*, \mu ) \ri ) = 1 - \f{\de}{2 s}$
for all $\ell \in \mscr{S}$.

In Case B, we have $\{ \mu \leq \ovl{X}_n, \; \mscr{M}_{\mrm{B}} (
\mu + \f{n}{n \vee m_\ell} ( \ovl{X}_n - \mu ), \mu ) \leq
\f{1}{m_\ell} \ln \f{\de}{2 s} \} = \{ \mu \leq \ovl{X}_n, \; \ln
\mu \leq \mscr{M}_{\mrm{B}} ( \mu + \f{n}{n \vee m_\ell} ( \ovl{X}_n
- \mu ), \mu ) \leq \f{1}{m_\ell} \ln \f{\de}{2 s} \} = \emptyset$.
It follows that $\{ \mu \leq L_n^\ell \} = \emptyset$ for all $n \in
\mscr{N}$. Therefore, $\Pr \{ L_n^\ell < \mu \; \tx{for all} \; n
\in \mscr{N} \} \geq 1 - \sum_{n \in \mscr{N}} \Pr \{ \mu \leq
L_n^\ell \} = 1$ for all $\ell \in \mscr{S}$, which implies that
$\Pr \{ L_n^\ell < \mu \; \tx{for all} \; n \in \mscr{N} \} = 1$ for
all $\ell \in \mscr{S}$. This completes the proof of the lemma.

\epf

\beL \la{lem3338b} For $n \in \mscr{N}$ and $\ell \in \mscr{S}$, define \[ U_n^\ell = \bec \sup \li \{ \nu \in (\ovl{X}_n, 1):
\mscr{M}_{\mrm{B}} \li ( \nu - \f{n}{n \vee m_\ell} \li ( \nu - \ovl{X}_n \ri ), \nu \ri ) > \f{1}{m_\ell} \ln \f{\de}{2 s} \ri \}  & \; \tx{for} \; \ovl{X}_n < 1,\\
1 & \; \tx{for} \; \ovl{X}_n = 1.  \eec \]  Then, $\Pr \{ U_n^\ell > \mu \; \tx{for all} \; n \in \mscr{N} \} \geq 1 - \f{\de}{2 s}$ for all
$\ell \in \mscr{S}$. \eeL

\bpf

First, we need to show that $U_n^\ell$ is well-defined.  Since $U_n^\ell = 1$ for $\ovl{X}_n = 1$, $U_n^\ell$ is well-defined provided that
$U_n^\ell$ exists for $0 \leq \ovl{X}_n < 1$.  Note that  $\lim_{\nu \downarrow y} \mscr{M}_{\mrm{B}} ( \nu - \f{n}{n \vee m_\ell} \li ( \nu - y
\ri ), \nu  ) = 0
> \f{1}{m_\ell} \ln \f{\de}{2 s}$ for $y \in [0, 1)$.  This fact together with Lemma \ref{lem2338b} imply the existence of $U_n^\ell$ for $0 \leq \ovl{X}_n
< 1$. So, $U_n^\ell$ is well-defined.  From the definition of $U_n^\ell$, it can be seen that \bee &  & \{ \mu \geq U_n^\ell, \; \ovl{X}_n = 1
\}  = \li \{ \mu \geq \ovl{X}_n,
\; \mscr{M}_{\mrm{B}} \li ( \mu - \f{n}{n \vee m_\ell} ( \mu - \ovl{X}_n ), \mu \ri ) \leq \f{1}{m_\ell} \ln \f{\de}{2 s}, \; \ovl{X}_n = 1 \ri \} = \emptyset,\\
&  & \{ \mu \geq U_n^\ell, \; 0 \leq \ovl{X}_n < 1 \} \subseteq  \li \{ \mu \geq \ovl{X}_n, \; \mscr{M}_{\mrm{B}} \li ( \mu - \f{n}{n \vee
m_\ell} ( \mu - \ovl{X}_n  ), \mu \ri ) \leq \f{1}{m_\ell} \ln \f{\de}{2 s}, \; 0 \leq \ovl{X}_n < 1 \ri \}. \eee This implies that $\{ \mu \geq
U_n^\ell \} \subseteq  \{ \mu \geq \ovl{X}_n, \; \mscr{M}_{\mrm{B}} ( \mu - \f{n}{n \vee m_\ell} ( \mu - \ovl{X}_n ), \mu ) \leq \f{1}{m_\ell}
\ln \f{\de}{2 s} \}$.

Next, consider $\Pr \{ U_n^\ell > \mu \; \tx{for all} \; n \in \mscr{N} \}$ for two cases as follows.

Case A: $(1- \mu)^{m_\ell} \leq \f{\de}{2 s}$.

Case B: $(1 - \mu)^{m_\ell} > \f{\de}{2 s}$.

In Case A, there must exist an $\vep^* \in (0, \mu]$ such that
$\mscr{M}_{\mrm{B}} \li ( \mu - \vep^*, \mu \ri ) = \f{1}{m_\ell}
\ln \f{\de}{2 s}$. Note that $\mscr{M}_{\mrm{B}} ( \mu - \ep, \mu )$
is decreasing with respect to $\ep \in (0, \mu)$. Therefore, from
the definitions of $U_n^\ell$ and $\vep^*$, we have that $\{ \mu
\geq U_n^\ell \} \subseteq  \{ \mu \geq \ovl{X}_n, \;
\mscr{M}_{\mrm{B}} ( \mu - \f{n}{n \vee m_\ell} ( \mu - \ovl{X}_n ),
\mu ) \leq \f{1}{m_\ell} \ln \f{\de}{2 s} \} \subseteq \{ \mu \geq
\ovl{X}_n, \; \f{n}{n \vee m_\ell} ( \mu - \ovl{X}_n ) \geq \vep^*
\} \subseteq \{ \ovl{X}_n \leq \mu - \f{(n \vee m_\ell) \vep^*}{n}
\}$. This implies that $\{ U_n^\ell > \mu \} \supseteq  \{ \ovl{X}_n
> \mu - \f{(n \vee m_\ell) \vep^*}{n}  \}$ for all $n \in \mscr{N}$.  Hence, $\{ U_n^\ell > \mu \; \tx{for all} \; n \in \mscr{N} \} \supseteq  \{ \ovl{X}_n >
\mu - \f{(n \vee m_\ell) \vep^*}{n} \; \tx{for all} \; n \in \mscr{N} \}$.  It follows from Lemma \ref{lem1338b} that $\Pr \{ U_n^\ell > \mu \;
\tx{for all} \; n \in \mscr{N} \} \geq \Pr \{  \ovl{X}_n
> \mu - \f{(n \vee m_\ell) \vep^*}{n} \; \tx{for all} \; n \in \mscr{N}  \} \geq 1 - \exp \li ( m_\ell \mscr{M}_{\mrm{B}} ( \mu - \vep^*, \mu ) \ri ) = 1 -
\f{\de}{2 s}$ for all $\ell \in \mscr{S}$.

In Case B, we have $\{ \mu \geq \ovl{X}_n, \; \mscr{M}_{\mrm{B}} (
\mu - \f{n}{n \vee m_\ell} ( \mu - \ovl{X}_n  ), \mu ) \leq
\f{1}{m_\ell} \ln \f{\de}{2 s} \} = \{ \mu \geq \ovl{X}_n, \; \ln (1
- \mu) \leq \mscr{M}_{\mrm{B}} ( \mu - \f{n}{n \vee m_\ell} ( \mu -
\ovl{X}_n  ), \mu ) \leq \f{1}{m_\ell} \ln \f{\de}{2 s} \} =
\emptyset$. It follows that $\{ \mu \geq U_n^\ell \} = \emptyset$
for all $n \in \mscr{N}$. Therefore, $\Pr \{ U_n^\ell
> \mu \; \tx{for all} \; n \in \mscr{N} \} \geq 1 - \sum_{n \in \mscr{N}} \Pr \{ \mu \geq U_n^\ell \} = 1$ for all
$\ell \in \mscr{S}$, which implies that $\Pr \{ U_n^\ell > \mu \; \tx{for all} \; n \in \mscr{N} \} = 1$ for all $\ell \in \mscr{S}$. This
completes the proof of the lemma.

\epf

\beL

\la{vip89}

{\small \bee &   &  \{ \ovl{X}_n - \vep  \leq  L_n^\ell \leq U_n^\ell \leq \ovl{X}_n + \vep \} \\
&  & = \li \{ \mscr{M}_{\mrm{B}} \li ( \ovl{X}_n + \vep - \f{n \vep}{n \vee m_\ell}, \ovl{X}_n + \vep \ri ) \leq \f{1}{m_\ell} \ln \f{\de}{2 s},
\; \mscr{M}_{\mrm{B}} \li ( \ovl{X}_n - \vep + \f{n \vep}{n \vee m_\ell}, \ovl{X}_n - \vep \ri ) \leq \f{1}{m_\ell} \ln \f{\de}{2 s} \ri \}
\eee} for all $n \in \mscr{N}$ and $\ell \in \mscr{S}$.

\eeL

\bpf  From the definitions of $\mscr{M}_{\mrm{B}}$ and $L_n^\ell$, it is clear that \be \la{coma} \li  \{ \mscr{M}_{\mrm{B}} \li ( \ovl{X}_n -
\vep + \f{n \vep}{n \vee m_\ell}, \ovl{X}_n - \vep \ri ) \leq \f{1}{m_\ell} \ln \f{\de}{2 s}, \; \ovl{X}_n - \vep \leq 0 \ri \} =  \{ \ovl{X}_n
- \vep \leq 0, \; \ovl{X}_n - \vep \leq L_n^\ell \} \ee for $n \in \mscr{N}$.  By Lemma \ref{lem2338} and the definition of $L_n^\ell$, \be
\la{comb} \li \{ \mscr{M}_{\mrm{B}} \li ( \ovl{X}_n - \vep + \f{n \vep}{n \vee m_\ell}, \ovl{X}_n - \vep \ri ) \leq \f{1}{m_\ell} \ln \f{\de}{2
s}, \;  \ovl{X}_n - \vep > 0 \ri \} = \{ 0 < \ovl{X}_n - \vep \leq L_n^\ell \} \ee for $n \in \mscr{N}$. It follows from (\ref{coma}) and
(\ref{comb}) that \be \la{part889a} \{ \ovl{X}_n - \vep \leq L_n^\ell \} = \li \{ \mscr{M}_{\mrm{B}} \li ( \ovl{X}_n - \vep + \f{n \vep}{n \vee
m_\ell}, \ovl{X}_n - \vep \ri ) \leq \f{1}{m_\ell} \ln \f{\de}{2 s} \ri \} \ee  for $n \in \mscr{N}$ and $\ell \in \mscr{S}$.

On the other hand, from the definitions of $\mscr{M}_{\mrm{B}}$ and $U_n^\ell$, it is clear that \be \la{coma2} \li \{ \mscr{M}_{\mrm{B}} \li (
\ovl{X}_n + \vep - \f{n \vep}{n \vee m_\ell}, \ovl{X}_n + \vep \ri )
 \leq \f{1}{m_\ell} \ln \f{\de}{2 s}, \; \ovl{X}_n + \vep \geq 1 \ri \} =  \{ \ovl{X}_n + \vep \geq 1, \; \ovl{X}_n +
\vep \geq U_n^\ell \} \ee for $n \in \mscr{N}$.  By Lemma \ref{lem2338b} and the definition of $U_n^\ell$, \be \la{comb2} \li  \{
\mscr{M}_{\mrm{B}} \li ( \ovl{X}_n + \vep - \f{n \vep}{n \vee m_\ell}, \ovl{X}_n + \vep \ri ) \leq \f{1}{m_\ell} \ln \f{\de}{2 s}, \;  \ovl{X}_n
+ \vep < 1 \ri \} = \{ 1 > \ovl{X}_n + \vep \geq U_n^\ell \} \ee for $n \in \mscr{N}$. It follows from (\ref{coma2}) and (\ref{comb2})  that \be
\la{part889b}  \{ \ovl{X}_n + \vep \geq U_n^\ell \} = \li \{ \mscr{M}_{\mrm{B}} \li ( \ovl{X}_n + \vep - \f{n \vep}{n \vee m_\ell}, \ovl{X}_n +
\vep \ri ) \leq \f{1}{m_\ell} \ln \f{\de}{2 s} \ri \} \ee for $n \in \mscr{N}$ and $\ell \in \mscr{S}$.  Combining (\ref{part889a}) and
(\ref{part889b}) completes the proof of the lemma.

\epf

\beL  \la{vip898} Let $0 < \vep < \f{1}{2}$ and $0 < r \leq 1$.  Then, \bel &   & \mscr{M}_{\mrm{B}} \li ( y + \vep - r \vep, y + \vep \ri )
\leq \mscr{M}_{\mrm{B}} \li ( y - \vep + r
\vep, y - \vep \ri ) \qqu \tx{for} \;  y \in \li [ \f{1}{2}, 1 \ri ], \la{showa}\\
&  &  \mscr{M}_{\mrm{B}} \li ( y + \vep - r \vep, y + \vep \ri ) \geq \mscr{M}_{\mrm{B}} \li ( y - \vep + r \vep, y - \vep \ri ) \qqu \tx{for}
\; y \in \li [ 0,  \f{1}{2} \ri ]. \la{showb} \eel \eeL

\bpf

First, we need to show (\ref{showa}). Note that $0 < \vep < \f{1}{2} \leq y \leq 1$.  For $y \geq 1 - \vep$, we have $0 < y - \vep < y - \vep +
r \vep \leq 1$ and $\mscr{M}_{\mrm{B}} \li ( y + \vep - r \vep, y + \vep \ri ) = - \iy < \mscr{M}_{\mrm{B}} \li ( y - \vep + r \vep, y - \vep
\ri )$.  Thus, to show (\ref{showa}), it suffices to show $\mscr{M}_{\mrm{B}} \li ( y + \vep - r \vep, y + \vep \ri ) \leq \mscr{M}_{\mrm{B}}
\li ( y - \vep + r \vep, y - \vep \ri )$ for $\f{1}{2} \leq y < 1 - \vep$.  Let $\se = y + \vep, \; z = y + \vep - r \vep$ and $\vse = y - \vep,
\; w = y - \vep + r \vep$.  A tedious computation shows that
 \bee &  & \f{ \pa  \mscr{M}_{\mrm{B}} \li ( y + \vep -
r \vep, y + \vep \ri )  }{\pa \vep} - \f{ \pa  \mscr{M}_{\mrm{B}} \li ( y - \vep + r \vep, y - \vep \ri )  }{\pa \vep} \\
&  & = ( 1 - r )  \ln \f{\se (1 - z)}{z (1 - \se)}  \f{\vse (1 - w)}{w (1 - \vse)} - r \vep  \li [ \f{1}{\se (1 - \se)} - \f{1}{
\vse (1 - \vse) } \ri ]\\
&  & \leq ( 1 - r ) \li [ \f{\se (1 - z)\vse (1 - w) }{z (1 - \se) w (1 - \vse)} - 1 \ri ] -  r \vep \li [ \f{1}{\se (1 - \se)} - \f{1}{ \vse (1
- \vse) } \ri ]\\
&  & = \f{ r^2 \vep^2 (2 y - 1) [ (r - 3) y^2 - (1 - r) \vep^2   ]   } {  (1 - \se) ( 1 - \vse) (y^2 - \vep^2) [ y^2 - (1 - r)^2 \vep^2 ]  }
\leq 0 \eee for $\f{1}{2} \leq y < 1 - \vep$.  Using this fact and the observation that $\mscr{M}_{\mrm{B}} \li ( y + \vep - r \vep, y + \vep
\ri ) \leq \mscr{M}_{\mrm{B}} \li ( y - \vep + r \vep, y - \vep \ri )$ for $\vep = 0$, we have that $\mscr{M}_{\mrm{B}} \li ( y + \vep - r \vep,
y + \vep \ri ) \leq \mscr{M}_{\mrm{B}} \li ( y - \vep + r \vep, y - \vep \ri )$ for $\f{1}{2} \leq y < 1 - \vep$.   Thus,  we have shown
(\ref{showa}).

Making use of (\ref{showa}) and  the relationship that $\mscr{M}_{\mrm{B}} (z, \se) = \mscr{M}_{\mrm{B}} (1 - z, 1 - \se)$ for $z \in [0, 1]$
and $\se \in (0, 1)$, we can conclude (\ref{showb}).

\epf

Making use of Lemmas \ref{vip89}, \ref{vip898} and the fact that $\mscr{M}_{\mrm{B}} (z, \se) = \mscr{M}_{\mrm{B}} (1 - z, 1 - \se)$ for $z \in
[0, 1]$ and $\se \in (0, 1)$, we have the following result as asserted by Lemma \ref{equ889}.

 \beL
 \la{equ889}

{\small \[ \{ \ovl{X}_n - \vep  \leq  L_n^\ell \leq U_n^\ell \leq \ovl{X}_n + \vep \} = \li \{ \mscr{M}_{\mrm{B}} \li ( \f{1}{2} - \li |
\f{1}{2} - \ovl{X}_n \ri | + \vep - \f{n \vep}{n \vee m_\ell}, \; \f{1}{2} - \li | \f{1}{2} - \ovl{X}_n \ri |  + \vep \ri ) \leq \f{1}{m_\ell}
\ln \f{\de}{2s}  \ri \}
\]} for all $n \in \mscr{N}$ and $\ell \in \mscr{S}$.

\eeL

\beL \la{confidenceabs}
Define $\bs{\mcal{L}}_n = \max_{\ell \in \mscr{S}} L_n^\ell$ and $\bs{\mcal{U}}_n = \min_{\ell \in \mscr{S}} U_n^\ell$
for $n \in \mscr{N}$.  Then, $\Pr \{  \bs{\mcal{L}}_n < \mu < \bs{\mcal{U}}_n  \; \tx{for all} \; n \in \mscr{N} \} \geq 1 - \de$.

\eeL

\bpf  Recall that $\Pr \{   L_n^\ell < \mu \; \tx{for all} \;  n \in \mscr{N} \} \geq 1 - \f{\de}{2 s}$ for $\ell \in \mscr{S}$.  By
Bonferroni's inequality, we have that $\Pr \{   L_n^\ell < \mu \; \tx{for all} \;  n \in \mscr{N} \; \tx{and all} \; \ell \in \mscr{S} \} \geq 1
- \f{\de}{2}$,  which implies that \be \la{parta} \Pr \{  \bs{\mcal{L}}_n < \mu \; \tx{for all} \; n \in \mscr{N} \} \geq 1 - \f{\de}{2}. \ee

On the other hand, note that $\Pr \{   U_n^\ell > \mu \; \tx{for all} \;  n \in \mscr{N} \} \geq 1 - \f{\de}{2 s}$ for $\ell \in \mscr{S}$. By
Bonferroni's inequality, we have that $\Pr \{   U_n^\ell > \mu \; \tx{for all} \;  n \in \mscr{N} \; \tx{and all} \; \ell \in \mscr{S} \} \geq 1
- \f{\de}{2}$,  which implies that \be \la{partb} \Pr \{  \bs{\mcal{U}}_n > \mu \; \tx{for all} \; n \in \mscr{N} \} \geq 1 - \f{\de}{2}. \ee
Combining (\ref{parta}) and (\ref{partb}) proves the lemma.
 \epf

\beL

\la{Finitestopabs} $\Pr \{ \ovl{X}_n - \vep \leq \bs{\mcal{L}}_n \leq \bs{\mcal{U}}_n \leq \ovl{X}_n + \vep \; \tx{for some} \;  n \in \mscr{N}
\} = 1$.

\eeL

\bpf

By Lemma \ref{equ889} and the definition of $\bs{\mcal{L}}_n$ and $\bs{\mcal{U}}_n$, we have \bee &  & \Pr \{ \ovl{X}_n - \vep \leq
\bs{\mcal{L}}_n \leq \bs{\mcal{U}}_n \leq \ovl{X}_n + \vep \; \tx{for some} \;  n \in \mscr{N} \}\\
&  & \geq \Pr \{ \ovl{X}_{m_s} - \vep \leq L_{m_s}^s \leq U_{m_s}^s \leq \ovl{X}_{m_s} + \vep \}\\
&  &  =  \Pr \li \{ \mscr{M}_{\mrm{B}} \li ( \f{1}{2} - \li | \f{1}{2} - \ovl{X}_n \ri |, \; \f{1}{2} - \li | \f{1}{2} - \ovl{X}_n \ri |  + \vep
\ri ) \leq \f{1}{m_\ell} \ln \f{\de}{2s} \ri \}. \eee

As a consequence of  the assumption that $m_s \geq \f{ \ln \f{2 s}{\de} }{ 2 \vep^2 }$, we have $\Pr \{  \mscr{M}_{\mrm{B}}  ( \f{1}{2} - |
\f{1}{2} - \ovl{X}_n  |, \; \f{1}{2} -  | \f{1}{2} - \ovl{X}_n  |  + \vep  ) \leq \f{1}{m_\ell} \ln \f{\de}{2s} \} = 1$, from which the lemma
follows.

 \epf

\bsk

Now we are in a position to prove that stopping rule A ensures the desired level of coverage probability.  From Lemma \ref{equ889} , we know
that the stopping rule is equivalent to ``continue sampling until $\{ \ovl{X}_n - \vep \leq L_n^\ell \leq U_n^\ell \leq \ovl{X}_n + \vep \}$ for
some $\ell \in \mscr{S}$ and $n \in \mscr{N}$''.  We claim that this stopping rule implies that ``continue sampling until $\{ \ovl{X}_n - \vep
\leq \bs{\mcal{L}}_n \leq \bs{\mcal{U}}_n \leq \ovl{X}_n + \vep \}$ for some $n \in \mscr{N}$''. To show this claim, we need to show
\[
\bigcup_{\ell \in \mscr{S}} \bigcup_{n \in \mscr{N}} \{ \ovl{X}_n - \vep  \leq  L_n^\ell \leq U_n^\ell \leq \ovl{X}_n + \vep \} \subseteq
\bigcup_{n \in \mscr{N}} \{ \ovl{X}_n - \vep \leq \bs{\mcal{L}}_n \leq \bs{\mcal{U}}_n \leq \ovl{X}_n + \vep \},
\]
which follows from the fact that {\small $\bigcup_{\ell \in \mscr{S}}  \{ \ovl{X}_n - \vep  \leq  L_n^\ell \leq U_n^\ell \leq \ovl{X}_n + \vep
\} \subseteq \{ \ovl{X}_n - \vep  \leq \bs{\mcal{L}}_n \leq \bs{\mcal{U}}_n \leq \ovl{X}_n + \vep \}$}  for every $n \in \mscr{N}$.  From Lemma
\ref{Finitestopabs}, we know that the sampling process will terminate at or before the $s$-th stage.   It follows from Lemma \ref{confidenceabs}
and Theorem \ref{General Inclusion Principle} that $\Pr \{ | \ovl{X}_{\mbf{n}} - \mu | < \vep \} \geq 1 - \de$.

\subsection{Derivation of Stopping Rule B}  \la{De_app_B}

Define  function $\mcal{M}_{\mrm{B}}(., .)$ such that
\[
\mcal{M}_{\mrm{B}} (z, \se) = \bec \f{ 9 (z - \se)^2 } {2 \li ( z + 2 \se \ri ) \li ( z + 2 \se - 3 \ri ) } & \; \tx{for} \; z \in [0, 1], \;
\se \in (0, 1),\\
- \iy & \; \tx{for} \; z \in (- \iy, \iy), \; \se \notin (0, 1). \eec
\]
We have established the following result.   \beL \la{equa998}

\bee &  & \li \{  \li [ \li | \ovl{X}_n  - \f{1}{2} \ri | - \vep + \f{n \vep }{3 ( n \vee m_\ell ) } \ri ]^2 \geq \f{1}{4} - \li ( \f{n}{n \vee
m_\ell} \ri )^2 \f{ m_\ell \vep^2} { 2 \ln \f{2 s}{\de} } \ri \}\\
&  &  = \li \{ \mcal{M}_{\mrm{B}} \li ( \f{1}{2} - \li | \f{1}{2} - \ovl{X}_n \ri | + \vep - \f{n \vep }{n \vee m_\ell}, \; \f{1}{2} - \li |
\f{1}{2} - \ovl{X}_n \ri | + \vep \ri ) \leq \f{1}{m_\ell} \ln \f{\de}{2 s} \ri \} \eee
 for all $n \in \mscr{N}$ and $\ell \in \{1, \cd, s \}$.

\eeL

\bpf Let $y \in [0, 1], \; r \in (0, 1], \; \al \in (0, 1)$ and $m \in \mscr{N}$.  For simplicity of notations, let $\se = \f{1}{2} - \li | y -
\f{1}{2} \ri | + \vep, \; z = \se - r \vep$ and $w = \f{z + 2 \se}{3}$.  Then, $\se - z = r \vep, \; w = \se - \f{r \vep}{3}$ and $w ( 1 - w) =
(\f{1}{2} - w + \f{1}{2}) [\f{1}{2} - ( \f{1}{2} - w ) ] = \f{1}{4} - ( \f{1}{2} - w )^2 > 0$. Moreover, $\se \in (0, 1), \; z \in [0, 1]$. It
follows that
\[
\mcal{M}_{\mrm{B}} (z, \se) \leq \f{\ln \al}{m}  \qu \Longleftrightarrow  \qu  \f{ (r \vep)^2 }{ \f{1}{4} - ( \f{1}{2} - w )^2   } \geq \f{2}{m}
\ln \f{1}{\al} \qu \Longleftrightarrow  \qu  \li ( \f{1}{2} - w \ri )^2  \geq \f{1}{4}  + \f{m (r \vep)^2} { 2 \ln \al},
\]
which implies that \bee &  &  \li [ \li | y  - \f{1}{2} \ri | - \vep + \f{r \vep }{3 } \ri ]^2 \geq \f{1}{4} + \f{ m (r \vep)^2} { 2 \ln \al }
\qu \Longleftrightarrow  \qu \mcal{M}_{\mrm{B}} \li ( \f{1}{2} - \li | \f{1}{2} - y \ri | + \vep - r \vep, \; \f{1}{2} - \li | \f{1}{2} - y \ri
| + \vep \ri ) \leq \f{\ln \al}{m}.  \eee This proves the lemma.

\epf

As a consequence of Lemma \ref{equa998}, stopping rule B is equivalent to  the following stopping rule:

Continue sampling until $\mcal{M}_{\mrm{B}} \li ( \f{1}{2} - \li | \f{1}{2} - \ovl{X}_n \ri | + \vep - \f{n \vep }{n \vee m_\ell}, \; \f{1}{2} -
\li | \f{1}{2} - \ovl{X}_n \ri | + \vep \ri ) \leq \f{1}{m_\ell} \ln \f{\de}{2 s}$ for some integers $n \in \mscr{N}$ and $\ell \in \{1, \cd, s
\}$.

As a consequence of Massart's inequality, \bee &  &  \li \{ \mcal{M}_{\mrm{B}} \li ( \f{1}{2} - \li | \f{1}{2} - \ovl{X}_n \ri | + \vep - \f{n
\vep }{n \vee m_\ell}, \; \f{1}{2} - \li | \f{1}{2} - \ovl{X}_n \ri | + \vep \ri ) \leq \f{1}{m_\ell} \ln \f{\de}{2 s} \ri \}\\
&  & \subseteq \li \{ \mscr{M}_{\mrm{B}} \li ( \f{1}{2} - \li | \f{1}{2} - \ovl{X}_n \ri | + \vep - \f{n \vep }{n \vee m_\ell}, \; \f{1}{2} -
\li | \f{1}{2} - \ovl{X}_n \ri | + \vep \ri ) \leq \f{1}{m_\ell} \ln \f{\de}{2 s} \ri \} \eee for all $n \in \mscr{N}$ and $\ell \in \{1, \cd, s
\}$.

Thus, by a similar method as that used in Appendix \ref{De_app_A} to  justify that stopping rule A guarantees  the desired level of coverage
probability, we can show that stopping rule B also ensures that $\Pr \{ | \ovl{X}_{\mbf{n}} - \mu | < \vep \} \geq 1 - \de$.

\subsection{Derivation of Stopping Rule C}  \la{De_app_C}

 We need some preliminary results.

\beL \la{comm9886}
\[
\mscr{M}_{\mrm{B}} \li ( \f{y}{1 + \vep} + \f{r \vep y}{1 + \vep}, \f{y}{1 + \vep} \ri )  \geq \mscr{M}_{\mrm{B}} \li ( \f{y}{1 - \vep} - \f{r
\vep y}{1 - \vep}, \f{y}{1 - \vep} \ri )
\]
for $y, \; r  \in (0, 1]$ and $\vep \in (0, 1)$.

\eeL

\bpf  In the case of $y \geq 1 - \vep$, we have $\mscr{M}_{\mrm{B}}
( \f{y}{1 - \vep} - \f{r \vep y}{1 - \vep}, \f{y}{1 - \vep} ) = -
\iy < \mscr{M}_{\mrm{B}} ( \f{y}{1 + \vep} + \f{r \vep y}{1 + \vep},
\f{y}{1 + \vep} )$.  Therefore, it suffices to show the lemma for
the case that $0 < y < 1 - \vep$.  For simplicity of notations,
define $H (\vep) = \mscr{M}_{\mrm{B}} \li ( \f{y ( 1 + r \vep)}{1 +
\vep}, \f{y}{1 + \vep} \ri )$.  Then,
\[
H (\vep) = \li [ 1 - \f{y ( 1 + r \vep)}{1 + \vep} \ri ] \ln \f{ 1 +
\f{\vep}{1 - y}  }{  1 + \f{ 1 - r y}{1 - y} \vep } - \f{y ( 1 + r
\vep)}{1 + \vep} \ln (1 + r \vep).
\]
By virtue of Taylor expansion series \bee \ln (1 + r \vep)  = \sum_{\ell = 1}^\iy (-1)^{\ell + 1} \f{r^\ell }{\ell}  \vep^\ell, \qqu  \ln \f{ 1
+ \f{\vep}{1 - y} }{ 1 + \f{ 1 - r y}{1 - y} \vep }  = \sum_{\ell = 1}^\iy (-1)^{\ell + 1} \f{1}{\ell} \li [ \f{1}{(1 - y)^\ell} - \li ( \f{ 1 -
r y}{1 - y} \ri )^\ell \ri ] \vep^\ell  \eee and a lengthy computation, we have \[ (1 - \vep^2) [ H(\vep) - H(- \vep) ] = 2 \sum_{k = 1}^\iy
C(r, y, 2 k + 1) \; (1 - y)^{- 2 k} \vep^{2 k + 1} \] where \bee C(r, y, \ell) & = & \f{r^{\ell - 2} (1 - y)^{\ell-1}}{\ell - 2} -  \f{r^\ell (1
- y)^{\ell-1}}{\ell} - \f{(1 - r y) (1 - y)}{\ell - 2} \li [ r^{\ell - 2} (1 - y)^{\ell-2} + 1 - \li ( 1 - r y
\ri )^{\ell - 2} \ri ] \\
&   &  + \f{1}{\ell} \li [ r^\ell (1 - y)^\ell + 1 - \li ( 1 - r y
\ri )^\ell \ri ] -  \f{(1 - r) y}{\ell - 1} \li [ r^{\ell - 1} (1 -
y)^{\ell-1} + 1 - \li ( 1 - r y \ri )^{\ell - 1} \ri ] \eee for
$\ell > 2$.  A tedious computation shows that {\small \bee \f{\pa
C(r, y, \ell)}{\pa r}  = \f{y}{\ell - 1} \li [ r^{\ell - 1} (1 -
y)^{\ell-1} + 1 - \li ( 1 - r y \ri )^{\ell - 1} \ri ] + \f{r y (1 -
y)}{\ell - 2} \li [ r^{\ell - 2} (1 - y)^{\ell-2} + 1 - \li ( 1 - r
y \ri )^{\ell - 2} \ri ] \geq 0. \eee} But $C(r, y, \ell) = 0$ for
$r = 0$.  Thus, $C(r, y, \ell) \geq 0$ for all $r \in [0, 1]$ and $y
\in [0, 1]$. This proves that $H(\vep) - H(- \vep) \geq 0$ for $y,
\; r  \in (0, 1]$ and $\vep \in (0, 1)$.  Thus the lemma is
established.

\epf

\beL

\la{vip89rev}

{\small \[ \li \{ \f{\ovl{X}_n}{1 + \vep}  \leq  L_n^\ell \leq U_n^\ell \leq \f{\ovl{X}_n}{1 - \vep} \ri \} \supseteq  \li \{ \ovl{X}_n > 0, \;
\mscr{M}_{\mrm{B}} \li ( \f{\ovl{X}_n}{1 + \vep} \li ( 1 + \f{n \vep}{n \vee m_\ell} \ri ), \f{\ovl{X}_n}{1 + \vep} \ri ) \leq \f{1}{m_\ell} \ln
\f{\de_\ell}{2} \ri \}
\]} for all $n \in \mscr{N}$ and $\ell \in \bb{N}$.

\eeL

\bpf  From the definitions of $\mscr{M}_{\mrm{B}}$ and $L_n^\ell$, it is clear that \be \la{comarev} \li \{ \mscr{M}_{\mrm{B}} \li (
\f{\ovl{X}_n}{1 + \vep} \li ( 1 + \f{n \vep}{n \vee m_\ell} \ri ), \f{\ovl{X}_n}{1 + \vep} \ri ) \leq \f{1}{m_\ell} \ln \f{\de_\ell}{2}, \;
\ovl{X}_n = 0 \ri \} = \li \{ \ovl{X}_n = 0, \; \f{\ovl{X}_n}{ 1 + \vep } \leq L_n^\ell \ri \} \ee for $n \in \mscr{N}$.  By Lemma \ref{lem2338}
and the definition of $L_n^\ell$, \be \la{combrev} \li \{ \mscr{M}_{\mrm{B}} \li ( \f{\ovl{X}_n}{1 + \vep} \li ( 1 + \f{n \vep}{n \vee m_\ell}
\ri ), \f{\ovl{X}_n}{1 + \vep} \ri ) \leq \f{1}{m_\ell} \ln \f{\de_\ell}{2}, \; \ovl{X}_n  > 0 \ri \} = \li \{ 0 < \f{\ovl{X}_n}{1 + \vep} \leq
L_n^\ell \ri \} \ee for $n \in \mscr{N}$.  It follows from (\ref{comarev}) and (\ref{combrev}) that {\small \be \la{good899a}
 \li \{
\f{\ovl{X}_n}{1 + \vep} \leq L_n^\ell \ri \} = \li \{ \mscr{M}_{\mrm{B}} \li ( \f{\ovl{X}_n}{1 + \vep} \li ( 1 + \f{n \vep}{n \vee m_\ell} \ri
), \f{\ovl{X}_n}{1 + \vep} \ri ) \leq \f{1}{m_\ell} \ln \f{\de_\ell}{2} \ri \} \ee} for all $n \in \mscr{N}$ and $\ell \in \bb{N}$.

From the definitions of $\mscr{M}_{\mrm{B}}$ and $U_n^\ell$, it is clear that \be \la{coma2rev} \li \{ \mscr{M}_{\mrm{B}} \li ( \f{\ovl{X}_n}{1
- \vep} \li ( 1 - \f{n \vep}{n \vee m_\ell} \ri ), \f{\ovl{X}_n}{1 - \vep} \ri ) \leq \f{1}{m_\ell} \ln \f{\de_\ell}{2}, \; \ovl{X}_n + \vep
\geq 1 \ri \} = \li \{ \ovl{X}_n + \vep \geq 1, \; \f{\ovl{X}_n}{ 1 - \vep } \geq U_n^\ell \ri \} \ee for $n \in \mscr{N}$.  By Lemma
\ref{lem2338b} and the definition of $U_n^\ell$, \be \la{comb2rev} \li \{ \mscr{M}_{\mrm{B}} \li ( \f{ \ovl{X}_n}{1 - \vep} \li ( 1 - \f{n
\vep}{n \vee m_\ell} \ri ), \f{\ovl{X}_n}{1 - \vep} \ri ) \leq \f{1}{m_\ell} \ln \f{\de_\ell}{2}, \; \vep < \ovl{X}_n + \vep < 1 \ri \} = \li \{
\vep < \ovl{X}_n + \vep < 1, \; \f{\ovl{X}_n}{ 1 - \vep } \geq U_n^\ell \ri \} \ee for $n \in \mscr{N}$.  It follows from (\ref{coma2rev}) and
(\ref{comb2rev}) that \be \la{good899b} \li \{ \ovl{X}_n > 0, \; \mscr{M}_{\mrm{B}} \li ( \f{ \ovl{X}_n}{1 - \vep} \li ( 1 - \f{n \vep}{n \vee
m_\ell} \ri ), \f{\ovl{X}_n}{1 - \vep} \ri ) \leq \f{1}{m_\ell} \ln \f{\de_\ell}{2} \ri \} \subseteq \li \{ \f{\ovl{X}_n}{ 1 - \vep } \geq
U_n^\ell \ri \} \ee for $n \in \mscr{N}$ and $\ell \in \bb{N}$. Finally, the proof of the lemma can be completed by combining (\ref{good899a}),
(\ref{good899b}) and using Lemma \ref{comm9886}.

\epf

\beL \la{Bonferoonrev} Define $\bs{\mcal{L}}_n = \sup_{\ell \in \bb{N}} L_n^\ell$ and $\bs{\mcal{U}}_n = \inf_{\ell \in \bb{N}} U_n^\ell$ for $n
\in \mscr{N}$. Then, $\Pr \{  \bs{\mcal{L}}_n \leq \mu \leq \bs{\mcal{U}}_n  \; \tx{for all} \; n \in \mscr{N} \} \geq 1 - \de$.

\eeL

\bpf Recall that $\Pr \{   L_n^\ell < \mu \; \tx{for all} \;  n \in \mscr{N} \} \geq 1 - \f{\de_\ell}{2}$ for $\ell \in \bb{N}$.  By
Bonferroni's inequality, we have that
\[
\Pr \{   L_n^\ell < \mu \; \tx{for all} \;  n \in \mscr{N} \; \tx{and} \; \ell = 1, \cd, k \} \geq 1 - \f{\sum_{\ell = 1}^k \de_\ell}{2} \]
 for
any $k \in \bb{N}$. By the continuity of the probability measure, we have \bee  \Pr \{   L_n^\ell < \mu \; \tx{for all} \;  n \in \mscr{N} \;
\tx{and} \; \ell \in \bb{N} \} & = & \lim_{k \to \iy} \Pr \{   L_n^\ell < \mu \; \tx{for all} \; n \in \mscr{N} \; \tx{and} \; \ell = 1, \cd, k
\}\\
& \geq & 1 - \lim_{k \to \iy} \f{\sum_{\ell = 1}^k \de_\ell}{2} = 1 - \f{\de}{2}, \eee which implies that \be \la{parta899} \Pr \{
\bs{\mcal{L}}_n \leq \mu \; \tx{for all} \; n \in \mscr{N} \} \geq 1 - \f{\de}{2}. \ee On the other hand, note that $\Pr \{   U_n^\ell > \mu \;
\tx{for all} \;  n \in \mscr{N} \} \geq 1 - \f{\de_\ell}{2}$ for $\ell \in \bb{N}$. By the continuity of the probability measure and
Bonferroni's inequality, we have that $\Pr \{ U_n^\ell
> \mu \; \tx{for all} \;  n \in \mscr{N} \; \tx{and all} \; \ell \in \bb{N} \} \geq 1 - \f{\de}{2}$,  which implies that \be \la{partb899} \Pr
\{  \bs{\mcal{U}}_n \geq \mu \; \tx{for all} \; n \in \mscr{N} \} \geq 1 - \f{\de}{2}. \ee Combining (\ref{parta899}) and (\ref{partb899})
proves the lemma.
 \epf

\beL \la{Finitestoprevbino} $\Pr \{ \f{\ovl{X}_n}{1 + \vep} \leq \bs{\mcal{L}}_n \leq \bs{\mcal{U}}_n \leq \f{\ovl{X}_n}{1 - \vep} \; \tx{for
some} \;  n \in \mscr{N} \} = 1$.

\eeL

\bpf

By the definition of $\bs{\mcal{L}}_n$ and $\bs{\mcal{U}}_n$, it is
sufficient to show that $\Pr \{ \f{\ovl{X}_{m_\ell} }{1 + \vep} \leq
L_{m_\ell}^\ell \leq U_{m_\ell}^\ell \leq \f{\ovl{X}_{m_\ell}}{1 -
\vep} \; \tx{for some} \;  \ell \in \bb{N}  \} = 1$.  From Lemma
\ref{vip89rev}, it can be seen that \bee &  &  \Pr \li \{
\f{\ovl{X}_{m_\ell} }{1 + \vep} \leq L_{m_\ell}^\ell \leq
U_{m_\ell}^\ell \leq \f{\ovl{X}_{m_\ell}}{1 - \vep} \; \tx{for some}
\;  \ell \in \bb{N} \ri \}\\
&   & \geq \Pr \li \{ \ovl{X}_{m_\ell} > 0, \; \mscr{M}_{\mrm{B}}
\li ( \ovl{X}_{m_\ell}, \f{\ovl{X}_{m_\ell}}{1 + \vep} \ri ) \leq
\f{1}{m_\ell} \ln \f{\de_\ell}{2}  \; \tx{for some} \; \ell \in
\bb{N} \ri \}. \eee  This inequality and Bonferroni's inequality
imply that \bee &  &  \Pr \li \{ \f{\ovl{X}_{m_\ell} }{1 + \vep}
\leq L_{m_\ell}^\ell \leq U_{m_\ell}^\ell \leq
\f{\ovl{X}_{m_\ell}}{1 - \vep} \; \tx{for some}
\;  \ell \in \bb{N} \ri \}\\
&   & \geq \lim_{\ell \to \iy} \Pr \{ \ovl{X}_{m_\ell}
> 0 \} + \lim_{\ell \to \iy} \Pr \li \{ \mscr{M}_{\mrm{B}}
\li ( \ovl{X}_{m_\ell}, \f{\ovl{X}_{m_\ell}}{1 + \vep} \ri ) \leq \f{1}{m_\ell} \ln \f{\de_\ell}{2} \ri \} - 1. \eee Since $\mu > 0$, it follows
from the law of large numbers that $\lim_{\ell \to \iy} \Pr \{ \ovl{X}_{m_\ell}
> 0 \} = 1$.   To complete the proof of the lemma, it remains  to show that
{\small $\lim_{\ell \to \iy} \Pr \li \{ \mscr{M}_{\mrm{B}} \li ( \ovl{X}_{m_\ell}, \f{\ovl{X}_{m_\ell}}{1 + \vep} \ri ) \leq \f{1}{m_\ell} \ln
\f{\de_\ell}{2} \ri \} = 1$}.  This is accomplished as follows.

Let $0 < \eta < 1$.  Noting that $\f{1}{m_\ell} \ln \f{\de_\ell}{2}$
is negative for any $\ell > 0$ and that {\small $\f{1}{m_\ell} \ln
\f{\de_\ell}{2} \to 0 > \mscr{M}_{\mrm{B}} ( \eta \mu, \f{\eta
\mu}{1 + \vep} )$} as $\ell \to \iy$, we have that there exists an
integer $\ka$ such that {\small $\mscr{M}_{\mrm{B}} ( \eta \mu,
\f{\eta \mu}{1 + \vep} ) < \f{1}{m_\ell} \ln \f{\de_\ell}{2}$} for
all $\ell \geq \ka$. For $\ell$ no less than such $\ka$, we claim
that $z < \eta \mu$ if {\small $\mscr{M}_{\mrm{B}} ( z, \f{z}{1 +
\vep} ) > \f{1}{m_\ell} \ln \f{\de_\ell}{2}$} and $z \in [0, 1]$. To
prove this claim, suppose, to get a contradiction, that $z \geq \eta
\mu$. Then, since {\small $\mscr{M}_{\mrm{B}} ( z, \f{z}{1 + \vep}
)$} is monotonically decreasing with respect to $z \in (0,1)$, we
have {\small $\mscr{M}_{\mrm{B}} (z, \f{z}{1 + \vep}  ) \leq
\mscr{M}_{\mrm{B}}  ( \eta \mu, \f{ \eta \mu}{1 + \vep}  ) <
\f{1}{m_\ell} \ln \f{\de_\ell}{2}$}, which is a contradiction.
Therefore, we have shown the claim and it follows that {\small $ \{
\mscr{M}_{\mrm{B}} ( \ovl{X}_{m_\ell}, \f{\ovl{X}_{m_\ell}}{1 +
\vep} ) > \f{1}{m_\ell} \ln \f{\de_\ell}{2} \} \subseteq \{
\ovl{X}_{m_\ell} < \eta \mu \}$} for $\ell \geq \ka$. So, {\small
\[  \Pr \li \{ \mscr{M}_{\mrm{B}}
\li ( \ovl{X}_{m_\ell}, \f{\ovl{X}_{m_\ell}}{1 + \vep} \ri )
> \f{1}{m_\ell} \ln \f{\de_\ell}{2} \ri \} \leq  \Pr \{ \ovl{X}_{m_\ell} < \eta
\mu  \} < \exp \li ( - \f{(1 - \eta)^2 \mu m_\ell}{2} \ri )
\]} for large enough $\ell$,  where the last inequality is due to the multiplicative Chernoff bound.  Since
$m_\ell \to \iy$ as $\ell \to \iy$, we have $\lim_{\ell \to \iy} \Pr \li \{ \mscr{M}_{\mrm{B}} \li ( \ovl{X}_{m_\ell}, \f{\ovl{X}_{m_\ell}}{1 +
\vep} \ri ) \leq \f{1}{m_\ell} \ln \f{\de_\ell}{2} \ri \} = 1$. This proves the lemma.

\epf

 \bsk

Now we are in a position to prove that stopping rule C ensures the desired level of coverage probability.  From Lemma \ref{vip89rev}, we know
that the stopping rule implies that ``continue sampling until $\{ \f{\ovl{X}_n}{1 + \vep}  \leq  L_n^\ell \leq U_n^\ell \leq \f{\ovl{X}_n}{1 -
\vep} \}$ for some $\ell \in \bb{N}$ and $n \in \mscr{N}$''.  We claim that  this stopping rule implies that ``continue sampling until $\{
\f{\ovl{X}_n}{1 + \vep} \leq \bs{\mcal{L}}_n \leq \bs{\mcal{U}}_n \leq \f{\ovl{X}_n}{1 - \vep} \}$ for some $n \in \mscr{N}$''. To show this
claim, we need to show
\[
\bigcup_{\ell \in \bb{N}} \bigcup_{n \in \mscr{N}} \li \{ \f{\ovl{X}_n}{ 1 + \vep }  \leq  L_n^\ell \leq U_n^\ell \leq \f{\ovl{X}_n}{1 - \vep}
\ri \} \subseteq \bigcup_{n \in \mscr{N}} \li \{ \f{\ovl{X}_n}{1 + \vep} \leq \bs{\mcal{L}}_n \leq \bs{\mcal{U}}_n \leq \f{\ovl{X}_n}{1 - \vep}
\ri \},
\]
which follows from the fact that {\small $\bigcup_{\ell \in \bb{N}}  \li \{ \f{\ovl{X}_n} {1 + \vep}  \leq  L_n^\ell \leq U_n^\ell \leq
\f{\ovl{X}_n}{1 - \vep} \ri \} \subseteq \li \{ \f{\ovl{X}_n}{ 1 + \vep} \leq \bs{\mcal{L}}_n \leq \bs{\mcal{U}}_n \leq \f{\ovl{X}_n}{1 - \vep}
\ri \}$}  for every $n \in \mscr{N}$.  From Lemma \ref{Finitestoprevbino}, we know that the sampling process will eventually terminate. It
follows from Lemma \ref{Bonferoonrev} and Theorem \ref{General Inclusion Principle} that $\Pr \{ | \ovl{X}_{\mbf{n}} - \mu | < \vep \mu \} \geq
1 - \de$.

\subsection{Derivation of Stopping Rule D}   \la{De_app_D}

We need some preliminary results.  As applications of Corollary 5 of \cite{Chenmax}, we have Lemmas \ref{lem1338geo} and \ref{lem1338bgeo}.

\beL \la{lem1338geo}

Let $\se \in (1, \iy)$.  Let $m \in \mscr{N}$ and $\vep > 0$.  Then,
\[
\Pr \li \{  \ovl{X}_n < \se + \f{(n \vee m) \vep}{n} \; \tx{for all} \; n \in \mscr{N} \ri \} \geq 1 - \exp \li ( m \mscr{M}_{\mrm{G}} ( \se +
\vep, \se ) \ri ).
\]
\eeL

\beL \la{lem1338bgeo}

Let $\se \in (1, \iy)$.  Let $m \in \mscr{N}$ and $\vep \in (0, \se)$.  Then,
\[
\Pr \li \{  \ovl{X}_n > \se - \f{(n \vee m) \vep}{n} \; \tx{for all} \; n \in \mscr{N} \ri \} \geq 1 - \exp \li ( m \mscr{M}_{\mrm{G}} ( \se -
\vep, \se ) \ri ).
\]
\eeL

\beL \la{lem2338geo} Let $y \geq 1$ and $0 < r \leq 1$.  Then, $\mscr{M}_{\mrm{G}}  ( \se + r  (  y - \se  ), \se  )$ increases with respect to
$\se \in (1, y)$. \eeL

\bpf  For simplicity of notations, let $z = \se + r \li (  y - \se
\ri )$.  It can checked that $\f{ \pa  \mscr{M}_{\mrm{G}} ( z, \se)
}{ \pa z } = \ln \f{z (1 - \se)}{\se (1 - z)}$ and $\f{ \pa
\mscr{M}_{\mrm{G}} ( z, \se) }{ \pa \se } = \f{\se - z}{\se (1 -
\se)}$.   By the chain rule of differentiation and the inequality
$\ln (1 + x) \leq x$ for $x > - 1$, we have  \bee \f{ \pa
\mscr{M}_{\mrm{G}} \li ( \se + r \li ( y - \se \ri ), \se \ri )}{\pa
\se} & = & \f{\se - z}{\se (1 - \se)} - ( 1 -
r) \ln \f{\se (1 - z)}{z (1 - \se)}\\
& \geq & \f{\se - z}{\se (1 - \se)}  -  \f{ ( 1 - r)  (\se - z)  }{ z (1 - \se) }  = \f{( \se - z ) [  z -  ( 1 - r) \se ]}{\se (1 - \se) z} =
\f{( \se - z ) r y}{\se (1 - \se) z}  > 0.  \eee  This proves the lemma.

\epf

\beL \la{lem2338bgeo} Let $y \geq 1$ and $0 < r \leq 1$.  Then, $\mscr{M}_{\mrm{G}} ( \se - r (  \se - y  ), \se )$ decreases with respect to
$\se \in (y, \iy)$.

\eeL

\bpf  For simplicity of notation, let $z = \se - r \li (  \se - y \ri )$.   It can checked that $\f{ \pa  \mscr{M}_{\mrm{G}} ( z, \se) }{ \pa z
} = \ln \f{z (1 - \se)}{\se (1 - z)}$ and $\f{ \pa \mscr{M}_{\mrm{G}} ( z, \se) }{ \pa \se } = \f{\se - z}{\se (1 - \se)}$.   Hence, \bee \f{
\pa \mscr{M}_{\mrm{G}} \li ( \se - r \li (  \se - y  \ri ), \se \ri )}{\pa \se} & = & ( 1 - r)  \ln \f{z (1 - \se)}{\se (1 - z)} - \f{z -
\se}{\se (1 - \se)} \leq  \f{( 1 - r)  (z - \se) }{ \se (1 - z) } - \f{ z
- \se  } {\se (1 - \se)}\\
&  &  = \f{( \se - z ) [ ( 1 - r) (\se - 1) -  ( z - 1)  ] } {\se (\se - 1) (z - 1)}  =  \f{( \se - z ) r ( 1- y) } {\se (\se - 1) (z - 1)} \leq
0.  \eee  This proves the lemma.

\epf

\beL \la{lem3338geo}  For $n \in \mscr{N}$ and $\ell \in \mscr{S}$, define \[ L_n^\ell = \bec \inf \li \{ \nu \in (1, \ovl{X}_n):
\mscr{M}_{\mrm{G}} \li ( \nu + \f{n}{n \vee m_\ell}  \li (
\ovl{X}_n - \nu \ri ), \nu \ri ) > \f{1}{m_\ell} \ln \f{\de}{2 s} \ri \}  & \; \tx{for} \; \ovl{X}_n > 1,\\
1 & \; \tx{for} \; \ovl{X}_n = 1.  \eec \]  Then, $\Pr \{ L_n^\ell < \se \; \tx{for all} \; n \in \mscr{N} \} \geq 1 - \f{\de}{2 s}$ for $\ell
\in \mscr{S}$. \eeL

\bpf

First, we need to show that $L_n^\ell$ is well-defined.  Since $L_n^\ell = 1$ for $\ovl{X}_n = 1$, $L_n^\ell$ is well-defined provided that
$L_n^\ell$ exists for $\ovl{X}_n > 1$.  Note that  $\lim_{\nu \uparrow y} \mscr{M}_{\mrm{G}} ( \nu + \f{n}{n \vee m_\ell}  \li ( y - \nu \ri ),
\nu  ) = 0 > \f{1}{m_\ell} \ln \f{\de}{2 s}$ for $y \in (1, \iy)$.  This fact together with Lemma \ref{lem2338geo} imply the existence of
$L_n^\ell$ for $\ovl{X}_n
> 1$. So, $L_n^\ell$ is well-defined.  From the definition of $L_n^\ell$, it can be seen that \bee &  & \{ \se \leq L_n^\ell, \; \ovl{X}_n = 1 \}  = \li \{ \se
\leq \ovl{X}_n, \;
\mscr{M}_{\mrm{G}} \li ( \se + \f{n}{n \vee m_\ell}  ( \ovl{X}_n - \se ), \se \ri ) \leq \f{1}{m_\ell} \ln \f{\de}{2 s}, \; \ovl{X}_n = 1 \ri \} = \emptyset,\\
&  & \{ \se \leq L_n^\ell, \; \ovl{X}_n > 1 \} \subseteq  \li \{ \se \leq \ovl{X}_n, \; \mscr{M}_{\mrm{G}} \li ( \se + \f{n}{n \vee m_\ell}  (
\ovl{X}_n - \se ), \se \ri ) \leq \f{1}{m_\ell} \ln \f{\de}{2 s}, \; \ovl{X}_n > 1 \ri \}. \eee This implies that $\{ \se \leq L_n^\ell \}
\subseteq \{ \se \leq \ovl{X}_n, \; \mscr{M}_{\mrm{G}} ( \se + \f{n}{n \vee m_\ell}  ( \ovl{X}_n - \se ), \se ) \leq \f{1}{m_\ell} \ln \f{\de}{2
s} \}$.

Next, consider $\Pr \{ L_n^\ell < \se \; \tx{for all} \; n \in \mscr{N} \}$.  Since
\[
\lim_{t \to \iy} \li [  t \ln \f{t}{\se} - (t - 1) \ln \f{t - 1}{\se - 1}  \ri ]  = \lim_{t \to \iy} \li [  t \ln \f{t}{t - 1} - \ln (\se - 1) +
\ln (t - 1) - t \ln \f{\se}{\se - 1} \ri ] = - \iy,
\] there must exist an $\vep^* > 0$ such that
$\mscr{M}_{\mrm{G}} \li ( \se + \vep^*, \se \ri ) = \f{1}{m_\ell}
\ln \f{\de}{2 s}$. Note that $\mscr{M}_{\mrm{G}} ( \se + \ep, \se )$
is decreasing with respect to $\ep > 0$. Therefore, from the
definitions of $L_n^\ell$ and $\vep^*$, we have that $\{ \se \leq
L_n^\ell \} \subseteq  \{ \se \leq \ovl{X}_n, \; \mscr{M}_{\mrm{G}}
( \se + \f{n}{n \vee m_\ell}  ( \ovl{X}_n - \se ), \se ) \leq
\f{1}{m_\ell} \ln \f{\de}{2 s} \} \subseteq \{ \se \leq \ovl{X}_n,
\; \f{n}{n \vee m_\ell} ( \ovl{X}_n - \se ) \geq \vep^* \} \subseteq
\{ \ovl{X}_n \geq \se + \f{ (n \vee m_\ell)  \vep^* }{n} \}$. This
implies that $\{ L_n^\ell < \se \} \supseteq  \{ \ovl{X}_n < \se +
\f{ (n \vee m_\ell)  \vep^* }{n}  \}$ for all $n \in \mscr{N}$.
Hence, $\{ L_n^\ell < \se \; \tx{for all} \; n \in \mscr{N} \}
\supseteq \{ \ovl{X}_n < \se + \f{ (n \vee m_\ell)  \vep^* }{n} \;
\tx{for all} \; n \in \mscr{N} \}$.  It follows from Lemma
\ref{lem1338geo} that $\Pr \{ L_n^\ell < \se \; \tx{for all} \; n
\in \mscr{N} \} \geq \Pr \{ \ovl{X}_n < \se + \f{ (n \vee m_\ell)
\vep^* }{n} \; \tx{for all} \; n \in \mscr{N}  \} \geq 1 - \exp \li
( m_\ell \mscr{M}_{\mrm{G}} ( \se + \vep^*, \se ) \ri ) = 1 -
\f{\de}{2 s}$ for $\ell \in \mscr{S}$.  This completes the proof of
the lemma.

\epf

\beL \la{lem3338bgeo} For $n \in \mscr{N}$ and $\ell \in \mscr{S}$, define $U_n^\ell = \sup  \{ \nu \in (\ovl{X}_n, 1): \mscr{M}_{\mrm{G}} ( \nu
- \f{n}{n \vee m_\ell} ( \nu - \ovl{X}_n ), \nu ) > \f{1}{m_\ell} \ln \f{\de}{2 s} \}$.   Then, $\Pr \{ U_n^\ell > \se \; \tx{for all} \; n \in
\mscr{N} \} \geq 1 - \f{\de}{2 s}$ for $\ell \in \mscr{S}$.  \eeL

\bpf

First, we need to show that $U_n^\ell$ is well-defined.  Note that  $\lim_{\nu \downarrow y} \mscr{M}_{\mrm{G}}  ( \nu - \f{n}{n \vee m_\ell}
\li ( \nu - y \ri ), \nu  ) = 0 > \f{1}{m_\ell} \ln \f{\de}{2 s}$ for $y \in [1, \iy)$.  This fact together with Lemma \ref{lem2338bgeo} imply
the existence of $U_n^\ell$. So, $U_n^\ell$ is well-defined.  From the definition of $U_n^\ell$, it can be seen that $\{ \se \geq U_n^\ell \}
\subseteq  \{ \se \geq \ovl{X}_n, \; \mscr{M}_{\mrm{G}} ( \se - \f{n}{n \vee m_\ell} ( \se - \ovl{X}_n ), \se ) \leq \f{1}{m_\ell} \ln \f{\de}{2
s} \}$.

Next, consider $\Pr \{ U_n^\ell > \se \; \tx{for all} \; n \in \mscr{N} \}$ for two cases as follows.

Case A: $\se^{-m_\ell} \leq \f{\de}{2 s}$.

Case B: $\se^{-m_\ell} > \f{\de}{2 s}$.

In Case A, there must exist an $\vep^* \in (0, \se - 1]$ such that
$\mscr{M}_{\mrm{G}} \li ( \se - \vep^*, \se \ri ) = \f{1}{m_\ell}
\ln \f{\de}{2 s}$. Note that $\mscr{M}_{\mrm{G}} ( \se - \ep, \se )$
is decreasing with respect to $\ep \in (0, \se - 1)$. Therefore,
from the definitions of $U_n^\ell$ and $\vep^*$, we have that $\{
\se \geq U_n^\ell \} \subseteq  \{ \se \geq \ovl{X}_n, \;
\mscr{M}_{\mrm{G}} ( \se - \f{n}{n \vee m_\ell} ( \se - \ovl{X}_n ),
\se ) \leq \f{1}{m_\ell} \ln \f{\de}{2 s} \} \subseteq \{ \se \geq
\ovl{X}_n, \; \f{n}{n \vee m_\ell} ( \se - \ovl{X}_n ) \geq \vep^*
\} \subseteq \{ \ovl{X}_n \leq \se - \f{ (n \vee m_\ell)  \vep^*
}{n} \}$. This implies that $\{ U_n^\ell > \se \} \supseteq  \{
\ovl{X}_n
> \se - \f{ (n \vee m_\ell)  \vep^* }{n} \}$ for all $n \in \mscr{N}$.
Hence, $\{ U_n^\ell > \se \; \tx{for all} \; n \in \mscr{N} \}
\supseteq  \{ \ovl{X}_n > \se - \f{ (n \vee m_\ell)  \vep^* }{n} \;
\tx{for all} \; n \in \mscr{N} \}$.  It follows from Lemma
\ref{lem1338bgeo} that $\Pr \{ U_n^\ell > \se \; \tx{for all} \; n
\in \mscr{N} \} \geq \Pr \{ \ovl{X}_n > \se - \f{ (n \vee m_\ell)
\vep^* }{n} \; \tx{for all} \; n \in \mscr{N}  \} \geq 1 - \exp \li
( m_\ell \mscr{M}_{\mrm{G}} ( \se - \vep^*, \se ) \ri ) = 1 -
\f{\de}{2 s}$.

In Case B, we have $\{ \se \geq \ovl{X}_n, \; \mscr{M}_{\mrm{G}} (
\se - \f{n}{n \vee m_\ell} ( \se - \ovl{X}_n  ), \se ) \leq
\f{1}{m_\ell} \ln \f{\de}{2 s} \} = \{ \se \geq \ovl{X}_n, \; \ln
\f{1}{\se} \leq \mscr{M}_{\mrm{G}} ( \se - \f{n}{n \vee m_\ell} (
\se - \ovl{X}_n  ), \se ) \leq \f{1}{m_\ell} \ln \f{\de}{2 s} \} =
\emptyset$. It follows that $\{ \se \geq U_n^\ell \} = \emptyset$
for all $n \in \mscr{N}$. Therefore, $\Pr \{ U_n^\ell > \se \;
\tx{for all} \; n \in \mscr{N} \} \geq 1 - \sum_{n \in \mscr{N}} \Pr
\{ \se \geq U_n^\ell \} = 1$ for $\ell \in \mscr{S}$, which implies
that $\Pr \{ U_n^\ell > \se \; \tx{for all} \; n \in \mscr{N} \} =
1$ for $\ell \in \mscr{S}$. This completes the proof of the lemma.

\epf

\beL

\la{vip89geo}

{\small \bee &   &  \{ (1 - \vep) \ovl{X}_n \leq  L_n^\ell \leq U_n^\ell \leq (1 + \vep) \ovl{X}_n \} \\
&  & = \li \{ \mscr{M}_{\mrm{G}} \li ( \li ( 1 + \vep - \f{n \vep}{n \vee m_\ell} \ri ) \ovl{X}_n, (1 + \vep) \ovl{X}_n \ri ) \leq \f{\ln
\f{\de}{2 s}}{m_\ell}, \; \mscr{M}_{\mrm{G}} \li ( \li ( 1 - \vep + \f{n \vep}{n \vee m_\ell} \ri ) \ovl{X}_n, (1 - \vep) \ovl{X}_n \ri ) \leq
\f{\ln \f{\de}{2 s}}{m_\ell} \ri \} \eee} for all $n \in \mscr{N}$ and $\ell \in \mscr{S}$.

\eeL

\bpf  From the definitions of $\mscr{M}_{\mrm{G}}$ and $L_n^\ell$, it is clear that {\small \be \la{comageo}
 \li \{ \mscr{M}_{\mrm{G}} \li ( \li ( 1 - \vep + \f{n \vep}{n \vee m_\ell} \ri ) \ovl{X}_n, (1 - \vep) \ovl{X}_n \ri )
 \leq \f{1}{m_\ell} \ln \f{\de}{2 s}, \; (1 - \vep) \ovl{X}_n \leq 1 \ri \} =  \{ (1 - \vep) \ovl{X}_n \leq 1, \; (1 - \vep) \ovl{X}_n \leq L_n^\ell
\} \qqu \ee} for $n \in \mscr{N}$.  By Lemma \ref{lem2338geo} and the definition of $L_n^\ell$, \be \la{combgeo}
 \li \{ \mscr{M}_{\mrm{G}} \li ( \li ( 1 - \vep + \f{n \vep}{n \vee m_\ell} \ri ), (1 - \vep) \ovl{X}_n \ri )
 \leq \f{1}{m_\ell} \ln \f{\de}{2 s}, \;  (1 - \vep) \ovl{X}_n > 1 \ri \} = \{ 1 < (1 - \vep) \ovl{X}_n \leq L_n^\ell \} \ee for $n \in
\mscr{N}$.  It follows from (\ref{comageo}) and (\ref{combgeo}) that \be \la{OK89}
 \{ (1 - \vep) \ovl{X}_n \leq L_n^\ell \} = \li \{
\mscr{M}_{\mrm{G}} \li ( \li ( 1 - \vep + \f{n \vep}{n \vee m_\ell} \ri ) \ovl{X}_n, (1 - \vep) \ovl{X}_n \ri ) \leq \f{1}{m_\ell} \ln \f{\de}{2
s} \ri \} \ee for $n \in \mscr{N}$ and $\ell \in \mscr{S}$.   By Lemma \ref{lem2338bgeo} and the definition of $U_n^\ell$, \be \la{comb2geo} \li
\{ \mscr{M}_{\mrm{G}} \li ( \li ( 1 + \vep - \f{n \vep}{n \vee m_\ell} \ri ) \ovl{X}_n, (1 + \vep) \ovl{X}_n \ri ) \leq \f{1}{m_\ell} \ln
\f{\de}{2 s} \ri \} = \{ (1 + \vep) \ovl{X}_n \geq U_n^\ell \} \ee for $n \in \mscr{N}$ and $\ell \in \mscr{S}$.  Combining (\ref{OK89}) and
(\ref{comb2geo}) completes the proof of the lemma.

\epf

\beL
 \la{vip898geo}  $\mscr{M}_{\mrm{G}} \li ( (1 - \vep + r \vep) y, (1 - \vep) y \ri ) \leq \mscr{M}_{\mrm{G}} \li ( (1 + \vep - r \vep) y, (1 + \vep) y \ri )$ for
$\vep \in (0, 1), \; r \in (0, 1]$ and $y \geq 1$.

\eeL

\bpf  In the case of $1 \leq y \leq \f{1}{1 - \vep}$, we have $\mscr{M}_{\mrm{G}} ( (1 - \vep + r \vep) y, (1 - \vep) y ) = - \iy <
\mscr{M}_{\mrm{G}} ( (1 + \vep - r \vep) y, (1 + \vep) y  )$.  Therefore, it suffices to show the lemma for the case that $y > \f{1}{1 - \vep}$.
For simplicity of notations, let $\nu = (1 - \vep) y, \; \vse = (1 + \vep) y, \; z = (1 - \vep + r \vep) y$ and $w = (1 + \vep - r \vep) y$.
Note that $\mscr{M}_{\mrm{G}} \li ( z, \nu \ri ) = \mscr{M}_{\mrm{G}} \li ( w, \vse \ri )$ for $\vep = 0$ and
 \bee &  & \f{ \pa  \mscr{M}_{\mrm{G}} \li ( z, \nu \ri )  }{\pa \vep} - \f{ \pa  \mscr{M}_{\mrm{G}} \li ( w, \vse \ri )  }{\pa \vep} \\
&  & = y ( 1 - r)  \ln \f{\nu (1 - z)}{z (1 - \nu)} +  y \f{z - \nu}{\nu (1 - \nu)}
+ y ( 1 - r)  \ln \f{\vse (1 - w)}{w (1 - \vse)} + y \f{w - \vse}{\vse (1 - \vse)}\\
&  & = ( 1 - r) y \ln \f{\nu (1 - z) \vse (1 - w)}{z (1 - \nu) w (1 - \vse)} - r \vep y \li [ \f{1}{(1 + \vep) (1 - \vse)} -
\f{1}{(1 - \vep) (1 - \nu)} \ri ]\\
&  & \leq ( 1 - r) y \li [  \f{\nu (1 - z) \vse (1 - w)}{z (1 - \nu) w (1 - \vse)} - 1 \ri ] - r \vep y \li [ \f{1}{(1 + \vep) (1 -
\vse)} - \f{1}{(1 - \vep) (1 - \nu)} \ri ]\\
&  & = \f{ (r \vep)^2  [ r - 3 - (1 - r) \vep^2  ] y (2 y - 1) }{  (1 - \nu) ( 1 - \vse) (1 - \vep^2) [ 1 - (1 - r)^2 \vep^2 ] } \leq 0,
 \eee
where the last inequality is a consequence of  $r \in (0, 1]$ and $y
\geq 1$.  This proves the lemma.

  \epf

Making use of Lemmas \ref{vip89geo} and \ref{vip898geo}, we have the following result.

\beL
 \la{equ889geo}

{\small \[ \{ (1 - \vep) \ovl{X}_n  \leq  L_n^\ell \leq U_n^\ell \leq (1 + \vep) \ovl{X}_n \} = \li \{ \mscr{M}_{\mrm{G}} \li ( \li (1 + \vep -
\f{n \vep}{n \vee m_\ell} \ri ) \ovl{X}_n, (1 + \vep) \ovl{X}_n \ri ) \leq \f{1}{m_\ell} \ln \f{\de}{2 s} \ri \}
\]} for all $n \in \mscr{N}$ and $\ell \in \mscr{S}$.

\eeL

By a similar argument as that for proving Lemma \ref{confidenceinverse}, we have established the following result.

\beL \la{confidenceinverse} Define $\bs{\mcal{L}}_n = \max_{\ell \in \mscr{S}} L_n^\ell$ and $\bs{\mcal{U}}_n = \min_{\ell \in \mscr{S}}
U_n^\ell$ for $n \in \mscr{N}$.  Then, $\Pr \{  \bs{\mcal{L}}_n < \mu < \bs{\mcal{U}}_n  \; \tx{for all} \; n \in \mscr{N} \} \geq 1 - \de$.

\eeL

\beL \la{Finitestopgeo} $\Pr \{ (1 - \vep) \ovl{X}_n  \leq \bs{\mcal{L}}_n \leq \bs{\mcal{U}}_n \leq ( 1 + \vep ) \ovl{X}_n  \; \tx{for some} \;
n \in \mscr{N} \} = 1$.

\eeL

\bpf

By Lemma \ref{equ889geo} and the definition of $\bs{\mcal{L}}_n$ and $\bs{\mcal{U}}_n$, we have \bee &  & \Pr \{ (1 - \vep) \ovl{X}_n  \leq
\bs{\mcal{L}}_n \leq \bs{\mcal{U}}_n \leq ( 1 + \vep ) \ovl{X}_n  \; \tx{for some} \;  n \in \mscr{N} \}\\
&  & \geq \Pr \{ (1 - \vep) \ovl{X}_{m_s} \leq L_{m_s}^s \leq U_{m_s}^s \leq (1 + \vep) \ovl{X}_{m_s}  \} =  \Pr \li \{ \mscr{M}_{\mrm{G}} \li (
\ovl{X}_{m_s}, (1 + \vep) \ovl{X}_{m_s} \ri ) \leq \f{1}{m_s} \ln \f{\de}{2 s} \ri \}.  \eee

As a consequence of  the assumption that $m_s \geq \f{(1 + \vep) \ln \f{2 s}{\de} }{ (1 + \vep) \ln (1 + \vep) - \vep }$, we have $\Pr \{
\mscr{M}_{\mrm{G}}  ( \ovl{X}_{m_s}, (1 + \vep) \ovl{X}_{m_s} ) \leq \f{1}{m_s} \ln \f{\de}{2 s} \} = 1$, from which the lemma follows.

 \epf

Now we are in a position to prove that stopping rule D ensures the desired level of coverage probability.  From Lemma \ref{equ889geo} , we know
that the stopping rule is equivalent to ``continue sampling until $\{ (1 - \vep) \ovl{X}_n \leq L_n^\ell \leq U_n^\ell \leq (1 + \vep) \ovl{X}_n
\}$ for some $\ell \in \mscr{S}$ and $n \in \mscr{N}$''.  We claim that this stopping rule implies that ``continue sampling until $\{ (1 - \vep)
\ovl{X}_n \leq \bs{\mcal{L}}_n \leq (1 + \vep) \bs{\mcal{U}}_n \leq \ovl{X}_n \}$ for some $n \in \mscr{N}$''. To show this claim, we need to
show
\[
\bigcup_{\ell \in \mscr{S}} \bigcup_{n \in \mscr{N}} \{ (1 - \vep) \ovl{X}_n  \leq  L_n^\ell \leq U_n^\ell \leq (1 + \vep) \ovl{X}_n  \}
\subseteq \bigcup_{n \in \mscr{N}} \{ (1 - \vep) \ovl{X}_n \leq \bs{\mcal{L}}_n \leq  \bs{\mcal{U}}_n \leq (1 + \vep) \ovl{X}_n \},
\]
which follows from the fact that {\small $\bigcup_{\ell \in \mscr{S}}  \{ (1 - \vep) \ovl{X}_n  \leq  L_n^\ell \leq U_n^\ell \leq (1 + \vep)
\ovl{X}_n \} \subseteq \{ (1 - \vep) \ovl{X}_n \leq \bs{\mcal{L}}_n \leq \bs{\mcal{U}}_n \leq (1 + \vep) \ovl{X}_n \}$}  for every $n \in
\mscr{N}$.  From Lemma \ref{Finitestopgeo}, we know that the sampling process will terminate at or before the $s$-th stage. It follows from
Lemma \ref{confidenceinverse} and Theorem \ref{General Inclusion Principle} that $\Pr \{ (1 - \vep) \ovl{X}_{\mbf{n}} < \se < (1 + \vep)
\ovl{X}_{\mbf{n}}  \} \geq 1 - \de$.

\subsection{Derivation of Stopping Rule E}  \la{De_app_E}

We need some preliminary results.

 \beL
\la{lempos}
 Let $y \geq 0$ and $0 < r \leq 1$.  Then, $\mscr{M}_{\mrm{P}} ( \lm + r (  y - \lm ), \lm )$ increases with respect to $\lm < y$.
 Similarly,  $\mscr{M}_{\mrm{P}} ( \lm - r ( \lm -  y ), \lm )$ decreases
with respect to $\lm > y$.

\eeL

\bpf  Note that $\mscr{M}_{\mrm{P}} ( z, \lm) = z - \lm  + z \ln
\f{\lm}{z}$.  It can be checked that $\f{ \pa  \mscr{M}_{\mrm{P}} (
z, \lm) }{ \pa z } = \ln \f{\lm}{z}$ and $\f{ \pa \mscr{M}_{\mrm{P}}
( z, \lm) }{ \pa \lm } = \f{z}{\lm} - 1$.  For simplicity of
notations, let $u = \lm + r \li (  y - \lm \ri )$ and $v = \lm - r
\li ( \lm - y \ri )$.   By the chain rule of differentiation,
{\small \bee &  & \f{ \pa \mscr{M}_{\mrm{P}} \li ( \lm + r \li ( y -
\lm \ri ), \lm \ri )}{\pa \lm}  =  \f{u - \lm}{\lm} - ( 1 - r)
\ln \f{u}{\lm} \geq \f{u - \lm}{\lm} -  ( 1 - r) \f{u - \lm}{\lm} = \f{r (u - \lm) }{\lm} \geq 0, \\
&  & \f{ \pa \mscr{M}_{\mrm{P}} \li ( \lm - r \li ( \lm - y \ri ), \lm \ri )}{\pa \lm} = ( 1 - r)  \ln \f{\lm}{v} + \f{v - \lm}{\lm} \leq ( 1 -
r) \f{\lm - v}{v} - \f{\lm - v}{\lm} = \f{ (v - \lm) y}{v \lm} \leq 0. \eee} This proves the lemma.

\epf

\beL \la{vip898pos}

$\mscr{M}_{\mrm{P}} \li ( y + \vep - r \vep, y + \vep \ri ) \geq \mscr{M}_{\mrm{P}} \li ( y - \vep + r \vep, y - \vep \ri )$ for $\vep > 0, \; y
\geq 0$ and $r \in (0, 1]$.

\eeL

\bpf  In the case of $0 \leq y \leq \vep$, we have $\mscr{M}_{\mrm{P}} \li ( y - \vep + r \vep, y - \vep \ri ) = - \iy < \mscr{M}_{\mrm{P}} \li
( y + \vep - r \vep, y + \vep \ri )$.  Therefore, it suffices to show the lemma for the case that $y > \vep$.  For simplicity of notations, let
$\se = y + \vep, \; \vse = y - \vep, \;  z = y + \vep - r \vep$ and $w = y - \vep + r \vep$. Note that
 \bee &  & \f{ \pa  \mscr{M}_{\mrm{P}} \li ( z, \se \ri )  }{\pa \vep} - \f{ \pa  \mscr{M}_{\mrm{P}} \li ( w, \vse \ri )  }{\pa \vep} \\
&  & = r \vep \li ( \f{1}{\vse} - \f{1}{ \se} \ri ) - ( 1 - r)  \ln \f{z w}{\se \vse} \geq  r \vep \li ( \f{1}{\vse} - \f{1}{ \se} \ri ) - ( 1 -
r) \li (  \f{z  w}{\se \vse } - 1 \ri ) = \f{ (r \vep)^2 (3 - r) }{ \se \vse } \geq 0.  \eee The proof of the lemma can be completed by making
use of this result and the observation that $\mscr{M}_{\mrm{P}} ( z, \se ) = \mscr{M}_{\mrm{P}} ( w, \vse)$ for $\vep = 0$.

\epf

As applications of Corollary 5 of \cite{Chenmax}, we have Lemmas \ref{lem1338pos} and \ref{lem1338bpos}.

\beL \la{lem1338pos}

Let $\lm \in (0, \iy)$.  Let $m \in \mscr{N}$ and $\vep > 0$.  Then,
\[
\Pr \li \{  \ovl{X}_n < \lm + \f{(m \vee n) \vep} {n} \; \tx{for
all} \; n \in \mscr{N} \ri \} \geq 1 - \exp \li ( m
\mscr{M}_{\mrm{P}} ( \lm + \vep, \lm ) \ri ).
\]
\eeL

\beL \la{lem1338bpos}

Let $\lm \in (0, \iy)$.  Let $m \in \mscr{N}$ and $\vep \in (0, \lm)$.  Then,
\[
\Pr \li \{  \ovl{X}_n > \lm - \f{(m \vee n) \vep} {n} \; \tx{for
all} \; n \in \mscr{N} \ri \} \geq 1 - \exp \li ( m
\mscr{M}_{\mrm{P}} ( \lm - \vep, \lm ) \ri ).
\]
\eeL

\beL \la{lem3338pos}  For $n \in \mscr{N}$ and $\ell \in \bb{N}$, define \[ L_n^\ell = \bec \inf \li \{ \nu \in (0, \ovl{X}_n):
\mscr{M}_{\mrm{P}} \li ( \nu + \f{n}{n \vee m_\ell}  \li (
\ovl{X}_n - \nu \ri ), \nu \ri ) > \f{1}{m_\ell} \ln \f{\de_\ell}{2} \ri \}  & \; \tx{for} \; \ovl{X}_n > 0,\\
0 & \; \tx{for} \; \ovl{X}_n = 0.  \eec \]  Then, $\Pr \{ L_n^\ell < \lm \; \tx{for all} \; n \in \mscr{N} \} \geq 1 - \f{\de_\ell}{2}$ for
$\ell \in \bb{N}$. \eeL

\bpf

First, we need to show that $L_n^\ell$ is well-defined.  Since
$L_n^\ell = 0$ for $\ovl{X}_n = 0$, $L_n^\ell$ is well-defined
provided that $L_n^\ell$ exists for $\ovl{X}_n > 0$.  Note that
$\lim_{\nu \uparrow y} \mscr{M}_{\mrm{P}} ( \nu + \f{n}{n \vee
m_\ell}  \li ( y - \nu \ri ), \nu ) = 0 > \f{1}{m_\ell} \ln
\f{\de_\ell}{2}$ for $y > 0$.  This fact together with Lemma
\ref{lempos} imply the existence of $L_n^\ell$ for $\ovl{X}_n > 0$.
So, $L_n^\ell$ is well-defined.  From the definition of $L_n^\ell$,
it can be seen that \bee &  & \{ \lm \leq L_n^\ell, \; \ovl{X}_n = 0
\} = \li \{ \lm \leq \ovl{X}_n, \;
\mscr{M}_{\mrm{P}} \li ( \lm + \f{n}{n \vee m_\ell}  ( \ovl{X}_n - \lm ), \lm \ri ) \leq \f{1}{m_\ell} \ln \f{\de_\ell}{2}, \; \ovl{X}_n = 0 \ri \} = \emptyset,\\
&  & \{ \lm \leq L_n^\ell, \; \ovl{X}_n > 0 \} \subseteq  \li \{ \lm \leq \ovl{X}_n, \; \mscr{M}_{\mrm{P}} \li ( \lm + \f{n}{n \vee m_\ell}  (
\ovl{X}_n - \lm ), \lm \ri ) \leq \f{1}{m_\ell} \ln \f{\de_\ell}{2}, \; \ovl{X}_n > 0 \ri \}. \eee This implies that $\{ \lm \leq L_n^\ell \}
\subseteq  \{ \lm \leq \ovl{X}_n, \; \mscr{M}_{\mrm{P}} ( \lm + \f{n}{n \vee m_\ell}  ( \ovl{X}_n - \lm ), \lm ) \leq \f{1}{m_\ell} \ln
\f{\de_\ell}{2} \}$.

Next, consider $\Pr \{ L_n^\ell < \lm \; \tx{for all} \; n \in
\mscr{N} \}$.  Since $\lim_{t \to \iy} t (\ln \f{\lm}{t} - 1) = -
\iy$, there must exist an $\vep^* > 0$ such that $\mscr{M}_{\mrm{P}}
\li ( \lm + \vep^*, \lm \ri ) = \f{1}{m_\ell} \ln \f{\de_\ell}{2}$.
Note that $\mscr{M}_{\mrm{P}} ( \lm + \ep, \lm )$ is decreasing with
respect to $\ep > 0$. Therefore, from the definitions of $L_n^\ell$
and $\vep^*$, we have that $\{ \lm \leq L_n^\ell \} \subseteq  \{
\lm \leq \ovl{X}_n, \; \mscr{M}_{\mrm{P}} ( \lm + \f{n}{n \vee
m_\ell}  ( \ovl{X}_n - \lm ), \lm ) \leq \f{1}{m_\ell} \ln
\f{\de_\ell}{2} \} \subseteq \{ \lm \leq \ovl{X}_n, \; \f{n}{n \vee
m_\ell}  ( \ovl{X}_n - \lm ) \geq \vep^* \} \subseteq \{ \ovl{X}_n
\geq \lm + \f{ (n \vee m_\ell)  \vep^* }{n} \}$. This implies that
$\{ L_n^\ell < \lm \} \supseteq  \{ \ovl{X}_n < \lm + \f{ (n \vee
m_\ell)  \vep^* }{n} \}$ for all $n \in \mscr{N}$. Hence, $\{
L_n^\ell < \lm \; \tx{for all} \; n \in \mscr{N} \} \supseteq  \{
\ovl{X}_n < \lm + \f{ (n \vee m_\ell)  \vep^* }{n} \; \tx{for all}
\; n \in \mscr{N} \}$. It follows from Lemma \ref{lem1338pos} that
$\Pr \{ L_n^\ell < \lm \; \tx{for all} \; n \in \mscr{N} \} \geq \Pr
\{ \ovl{X}_n < \lm + \f{ (n \vee m_\ell)  \vep^* }{n} \; \tx{for
all} \; n \in \mscr{N}  \} \geq 1 - \exp \li ( m_\ell
\mscr{M}_{\mrm{P}} ( \lm + \vep^*, \lm ) \ri ) = 1 -
\f{\de_\ell}{2}$ for $\ell \in \bb{N}$. This completes the proof of
the lemma.

\epf

\beL \la{lem3338bpos} For $n \in \mscr{N}$ and $\ell \in \bb{N}$, define $U_n^\ell = \sup \{ \nu \in (\ovl{X}_n, \iy): \mscr{M}_{\mrm{P}} ( \nu
- \f{n}{n \vee m_\ell} ( \nu - \ovl{X}_n ), \nu ) > \f{1}{m_\ell} \ln \f{\de_\ell}{2} \}$.   Then, $\Pr \{ U_n^\ell > \lm \; \tx{for all} \; n
\in \mscr{N} \} \geq 1 - \f{\de_\ell}{2}$ for $\ell \in \bb{N}$.  \eeL

\bpf

First, we need to show that $U_n^\ell$ is well-defined.  Note that  $\lim_{\nu \downarrow y} \mscr{M}_{\mrm{P}} ( \nu - \f{n}{n \vee m_\ell} \li
( \nu - y \ri ), \nu  ) = 0 > \f{1}{m_\ell} \ln \f{\de_\ell}{2}$ for $y \in [0, \iy)$.  This fact together with Lemma \ref{lempos} imply the
existence of $U_n^\ell$. So, $U_n^\ell$ is well-defined.  From the definition of $U_n^\ell$, it can be seen that $\{ \lm \geq U_n^\ell \}
\subseteq  \{ \lm \geq \ovl{X}_n, \; \mscr{M}_{\mrm{P}} ( \lm - \f{n}{n \vee m_\ell} ( \lm - \ovl{X}_n ), \lm ) \leq \f{1}{m_\ell} \ln
\f{\de_\ell}{2} \}$.

Next, consider $\Pr \{ U_n^\ell > \lm \; \tx{for all} \; n \in \mscr{N} \}$ for two cases as follows.

Case A: $\exp(- m_\ell \lm) \leq \f{\de_\ell}{2}$.

Case B: $\exp(- m_\ell \lm) > \f{\de_\ell}{2}$.

In Case A, there must exist an $\vep^* \in (0, \lm ]$ such that
$\mscr{M}_{\mrm{P}} \li ( \lm - \vep^*, \lm \ri ) = \f{1}{m_\ell}
\ln \f{\de_\ell}{2}$. Note that $\mscr{M}_{\mrm{P}} ( \lm - \ep, \lm
)$ is decreasing with respect to $\ep \in (0, \lm)$. Therefore, from
the definitions of $U_n^\ell$ and $\vep^*$, we have that $\{ \lm
\geq U_n^\ell \} \subseteq  \{ \lm \geq \ovl{X}_n, \;
\mscr{M}_{\mrm{P}} ( \lm - \f{n}{n \vee m_\ell} ( \lm - \ovl{X}_n ),
\lm ) \leq \f{1}{m_\ell} \ln \f{\de_\ell}{2} \} \subseteq \{ \lm
\geq \ovl{X}_n, \; \f{n}{n \vee m_\ell} ( \lm - \ovl{X}_n ) \geq
\vep^* \} \subseteq \{ \ovl{X}_n \leq \lm - \f{ (n \vee m_\ell)
\vep^* }{n} \}$. This implies that $\{ U_n^\ell > \lm \} \supseteq
\{ \ovl{X}_n > \lm - \f{ (n \vee m_\ell)  \vep^* }{n}  \}$ for all
$n \in \mscr{N}$. Hence, $\{ U_n^\ell > \lm \; \tx{for all} \; n \in
\mscr{N} \} \supseteq  \{ \ovl{X}_n > \lm - \f{ (n \vee m_\ell)
\vep^* }{n} \; \tx{for all} \; n \in \mscr{N} \}$. It follows from
Lemma \ref{lem1338bpos} that $\Pr \{ U_n^\ell > \lm \; \tx{for all}
\; n \in \mscr{N} \} \geq \Pr \{ \ovl{X}_n
> \lm - \f{ (n \vee m_\ell)  \vep^* }{n} \; \tx{for all} \; n \in \mscr{N}  \}
\geq 1 - \exp \li ( m_\ell \mscr{M}_{\mrm{P}} ( \lm - \vep^*, \lm ) \ri ) =
1 - \f{\de_\ell}{2}$ for $\ell \in \bb{N}$.

In Case B, we have $\{ \lm \geq \ovl{X}_n, \; \mscr{M}_{\mrm{P}} (
\lm - \f{n}{n \vee m_\ell} ( \lm - \ovl{X}_n  ), \lm ) \leq
\f{1}{m_\ell} \ln \f{\de_\ell}{2} \} = \{ \lm \geq \ovl{X}_n, \; -
\lm \leq \mscr{M}_{\mrm{P}} ( \lm - \f{n}{n \vee m_\ell} ( \lm -
\ovl{X}_n  ), \lm ) \leq \f{1}{m_\ell} \ln \f{\de_\ell}{2} \} =
\emptyset$. It follows that $\{ \lm \geq U_n^\ell \} = \emptyset$
for all $n \in \mscr{N}$. Therefore, $\Pr \{ U_n^\ell > \lm \;
\tx{for all} \; n \in \mscr{N} \} \geq 1 - \sum_{n \in \mscr{N}} \Pr
\{ \lm \geq U_n^\ell \} = 1$, which implies that $\Pr \{ U_n^\ell >
\lm \; \tx{for all} \; n \in \mscr{N} \} = 1$ for $\ell \in \bb{N}$.
This completes the proof of the lemma.

\epf

\beL

\la{vip89pos}

{\small \bee &   &  \{ \ovl{X}_n - \vep  \leq  L_n^\ell \leq U_n^\ell \leq \ovl{X}_n + \vep \} \\
&  & = \li \{ \mscr{M}_{\mrm{P}} \li ( \ovl{X}_n + \vep - \f{n \vep}{n \vee m_\ell}, \ovl{X}_n + \vep \ri ) \leq \f{1}{m_\ell} \ln
\f{\de_\ell}{2}, \; \mscr{M}_{\mrm{P}} \li ( \ovl{X}_n - \vep + \f{n \vep}{n \vee m_\ell}, \ovl{X}_n - \vep \ri ) \leq \f{1}{m_\ell} \ln
\f{\de_\ell}{2} \ri \} \eee} for all $n \in \mscr{N}$ and $\ell \in \bb{N}$.

\eeL

\bpf  From the definitions of $\mscr{M}_{\mrm{P}}$ and $L_n^\ell$, it is clear that {\small \be \la{comapos} \li \{ \mscr{M}_{\mrm{P}} \li (
\ovl{X}_n - \vep + \f{n \vep}{n \vee m_\ell}, \ovl{X}_n - \vep \ri ) \leq \f{1}{m_\ell} \ln \f{\de_\ell}{2}, \; \ovl{X}_n - \vep \leq 0 \ri \} =
\{ \ovl{X}_n - \vep \leq 0, \; \ovl{X}_n - \vep \leq L_n^\ell \} \ee} for $n \in \mscr{N}$.  By Lemma \ref{lempos} and the definition of
$L_n^\ell$, \be \la{combpos} \li \{ \mscr{M}_{\mrm{P}} \li ( \ovl{X}_n - \vep + \f{n \vep}{n \vee m_\ell}, \ovl{X}_n - \vep \ri ) \leq
\f{1}{m_\ell} \ln \f{\de_\ell}{2}, \;  \ovl{X}_n - \vep > 0 \ri \} = \{ 0 < \ovl{X}_n - \vep \leq L_n^\ell \} \ee for $n \in \mscr{N}$. It
follows from (\ref{comapos}) and (\ref{combpos}) that \be \la{hao1} \{ \ovl{X}_n - \vep \leq L_n^\ell \} = \li \{ \mscr{M}_{\mrm{P}} \li (
\ovl{X}_n - \vep + \f{n \vep}{n \vee m_\ell}, \ovl{X}_n - \vep \ri ) \leq \f{1}{m_\ell} \ln \f{\de_\ell}{2} \ri \} \ee for $n \in \mscr{N}$.

By Lemma \ref{lempos} and the definition of $U_n^\ell$,  \be \la{hao2} \{ \ovl{X}_n + \vep \geq U_n^\ell \} = \li \{ \mscr{M}_{\mrm{P}} \li (
\ovl{X}_n + \vep - \f{n \vep}{n \vee m_\ell}, \ovl{X}_n + \vep \ri ) \leq \f{1}{m_\ell} \ln \f{\de_\ell}{2} \ri \} \ee for $n \in \mscr{N}$.
Finally, combining (\ref{hao1}) and (\ref{hao2}) proves the lemma.

\epf

Making use of Lemmas \ref{vip89pos} and \ref{vip898pos}, we have the following result.

\beL
 \la{equ889pos}

{\small \[ \{ \ovl{X}_n - \vep  \leq  L_n^\ell \leq U_n^\ell \leq \ovl{X}_n + \vep \} = \li \{ \mscr{M}_{\mrm{P}} \li ( \ovl{X}_n + \vep - \f{n
\vep}{n \vee m_\ell}, \ovl{X}_n + \vep \ri ) \leq \f{1}{m_\ell} \ln \f{\de_\ell}{2} \ri \}
\]} for all $n \in \mscr{N}$ and $\ell \in \bb{N}$.

\eeL

By a similar argument as that for proving Lemma \ref{Bonferoonrev}, we have established the following result.

\beL \la{BonferoonrevPos} Define $\bs{\mcal{L}}_n = \sup_{\ell \in \bb{N}} L_n^\ell$ and $\bs{\mcal{U}}_n = \inf_{\ell \in \bb{N}} U_n^\ell$ for
$n \in \mscr{N}$. Then, $\Pr \{  \bs{\mcal{L}}_n \leq \mu \leq \bs{\mcal{U}}_n  \; \tx{for all} \; n \in \mscr{N} \} \geq 1 - \de$.

\eeL

\beL \la{Finitestopabspos} $\Pr \{ \ovl{X}_n - \vep  \leq  \bs{\mcal{L}}_n \leq \bs{\mcal{U}}_n \leq \ovl{X}_n + \vep   \; \tx{for some} \;  n
\in \mscr{N} \} = 1$.

\eeL

\bpf

By the definition of $\bs{\mcal{L}}_n$ and $\bs{\mcal{U}}_n$, it is sufficient to show that $\Pr \{ \ovl{X}_{m_\ell} - \vep \leq L_{m_\ell}^\ell
\leq U_{n_\ell}^\ell \leq \ovl{X}_{m_\ell} + \vep \; \tx{for some} \;  \ell \in \bb{N}  \} = 1$. In view of Lemma \ref{equ889pos}, this is
equivalent to show that $\Pr \{ \mscr{M}_{\mrm{P}} \li ( \ovl{X}_{m_\ell}, \ovl{X}_{m_\ell} + \vep \ri ) \leq \f{1}{m_\ell} \ln \f{\de_\ell}{2}
\; \tx{for some} \; \ell \in \bb{N} \} = 1$.  Note that $\Pr \{ \mscr{M}_{\mrm{P}} \li ( \ovl{X}_{m_\ell}, \ovl{X}_{m_\ell} + \vep \ri ) \leq
\f{1}{m_\ell} \ln \f{\de_\ell}{2}  \; \tx{for some} \; \ell \in \bb{N} \} \geq \lim_{\ell \to \iy} \Pr \{ \mscr{M}_{\mrm{P}} \li (
\ovl{X}_{m_\ell}, \ovl{X}_{m_\ell} + \vep \ri ) \leq \f{1}{m_\ell} \ln \f{\de_\ell}{2} \}$.  To complete the proof of the lemma, it remains to
show that $\lim_{\ell \to \iy} \Pr \{ \mscr{M}_{\mrm{P}} \li ( \ovl{X}_{m_\ell}, \ovl{X}_{m_\ell} + \vep \ri ) \leq \f{1}{m_\ell} \ln
\f{\de_\ell}{2} \} = 1$, which is accomplished as follows.

Let $0 < \eta < 1$. Noting that  {\small $\f{1}{m_\ell} \ln
\f{\de_\ell}{2} \to 0
> \mscr{M}_{\mrm{P}}  ( \f{\lm}{\eta}, \f{\lm}{\eta} + \vep )$} as
$\ell \to \iy$, we have that there exists an integer $\ka$ such that
{\small $\mscr{M}_{\mrm{P}}  ( \f{\lm}{\eta}, \f{\lm}{\eta} + \vep )
< \f{1}{m_\ell} \ln \f{\de_\ell}{2}$} for all $\ell \geq \ka$. For
$\ell$ no less than such $\ka$, we claim that $z > \f{\lm}{\eta}$ if
{\small $\mscr{M}_{\mrm{P}} ( z, z + \vep )
> \f{1}{m_\ell} \ln \f{\de_\ell}{2}$} and $z \in [0, \iy)$. To prove
this claim, suppose, to get a contradiction, that $z \leq
\f{\lm}{\eta}$. Then, since {\small $\mscr{M}_{\mrm{P}} ( z,  z +
\vep )$} is monotonically increasing with respect to $z > 0$, we
have {\small $\mscr{M}_{\mrm{P}} (z, z + \vep  ) \leq
\mscr{M}_{\mrm{P}} ( \f{\lm}{\eta}, \f{\lm}{\eta} + \vep ) <
\f{1}{m_\ell} \ln \f{\de_\ell}{2}$}, which is a contradiction.
Therefore, we have shown the claim and it follows that {\small $ \{
\mscr{M}_{\mrm{P}} ( \ovl{X}_{m_\ell}, \ovl{X}_{m_\ell} + \vep ) >
\f{1}{m_\ell} \ln \f{\de_\ell}{2} \} \subseteq \{ \ovl{X}_{m_\ell}
> \f{\lm}{\eta} \}$} for $\ell \geq \ka$. So, {\small $\Pr \{ \mscr{M}_{\mrm{P}} \li ( \ovl{X}_{m_\ell}, \ovl{X}_{m_\ell} +
\vep \ri ) > \f{1}{m_\ell} \ln \f{\de_\ell}{2}  \} \leq \Pr  \{ \ovl{X}_{m_\ell} > \f{\lm}{\eta} \} < \exp \li ( - c m_\ell  \ri )$,} where $c =
- \mscr{M}_{\mrm{P}} ( \f{\lm}{\eta}, \lm )$ and the last inequality is due to Chernoff bounds \cite{Chernoff}. Since $m_\ell \to \iy$ as $\ell
\to \iy$, we have $\lim_{\ell \to \iy} \Pr \{ \mscr{M}_{\mrm{P}} ( \ovl{X}_{m_\ell}, \ovl{X}_{m_\ell} + \vep ) \leq \f{1}{m_\ell} \ln
\f{\de_\ell}{2} \} = 1$. This proves the lemma.

\epf

\bsk

Now we are in a position to prove that stopping rule E ensures the desired level of coverage probability. From Lemma \ref{equ889pos}, we know
that the stopping rule implies that ``continue sampling until $\{ \ovl{X}_n - \vep  \leq  L_n^\ell \leq U_n^\ell \leq \ovl{X}_n + \vep \}$ for
some $\ell \in \bb{N}$ and $n \in \mscr{N}$''.  We claim that  this stopping rule implies that ``continue sampling until $\{ \ovl{X}_n - \vep
\leq \bs{\mcal{L}}_n \leq \bs{\mcal{U}}_n \leq \ovl{X}_n + \vep \}$ for some $n \in \mscr{N}$''. To show this claim, we need to show
\[
\bigcup_{\ell \in \bb{N}} \bigcup_{n \in \mscr{N}} \li \{ \ovl{X}_n - \vep \leq  L_n^\ell \leq U_n^\ell \leq \ovl{X}_n + \vep \ri \} \subseteq
\bigcup_{n \in \mscr{N}} \li \{ \ovl{X}_n - \vep \leq \bs{\mcal{L}}_n \leq \bs{\mcal{U}}_n \leq \ovl{X}_n + \vep \ri \},
\]
which follows from the fact that {\small $\bigcup_{\ell \in \bb{N}}  \li \{ \ovl{X}_n - \vep \leq  L_n^\ell \leq U_n^\ell \leq \ovl{X}_n + \vep
\ri \} \subseteq \li \{ \ovl{X}_n - \vep \leq \bs{\mcal{L}}_n \leq \bs{\mcal{U}}_n \leq \ovl{X}_n + \vep \ri \}$}  for every $n \in \mscr{N}$.
From Lemma \ref{Finitestopabspos}, we know that the sampling process will eventually terminate. It follows from Lemma \ref{BonferoonrevPos} and
Theorem \ref{General Inclusion Principle} that $\Pr \{ | \ovl{X}_{\mbf{n}} - \lm | < \vep \} \geq 1 - \de$.

\subsection{Derivation of Stopping Rule F}  \la{De_app_F}

 We need some preliminary results.

 \beL
 \la{vip898posrev}

 $\mscr{M}_{\mrm{P}} \li ( \f{(1 + r \vep) y}{1 + \vep}, \f{y}{1 + \vep} \ri ) >
 \mscr{M}_{\mrm{P}} \li ( \f{(1 - \vep - r \vep) y}{1 + \vep}, \f{y}{1 - \vep} \ri )$ for $\vep \in (0, 1), \; r \in (0, 1]$ and $y > 0$.
  \eeL

\bpf

For simplicity of notations, let $\se = \f{y}{1 + \vep}, \; \vse = \f{y}{1 - \vep}, \;  z = \f{y}{1 + \vep} + \f{r \vep y}{1 + \vep}$ and
  $w = \f{y}{1 - \vep} - \f{r \vep y}{1 - \vep}$.   Note that
 \bee &  & \f{ \pa  \mscr{M}_{\mrm{P}} \li ( z, \se \ri )  }{\pa \vep} - \f{ \pa  \mscr{M}_{\mrm{P}} \li ( w, \vse \ri )  }{\pa \vep} \\
&  & = \f{y}{(1 + \vep)^2} ( 1 - r)  \ln \f{z}{\se} - \f{y}{(1 + \vep)^2} \f{z - \se}{\se}
 + \f{y}{(1 - \vep)^2} ( 1 - r)  \ln \f{w}{\vse}
- \f{y}{(1 - \vep)^2} \f{w - \vse}{\vse}\\
 &  & = \f{ y (1 - r) }{ r^2 (1 - \vep^2)^2} [ h(r \vep) + h(- r \vep) ],
 \eee
 where $h (t) = (r - t)^2 \ln (1 + t) - \f{t}{1 - r} ( r - t)^2$.
 Using $\ln (1 + t) = t - \f{t^2}{2} + \f{t^3}{3} - \f{t^4}{4} + \cd$ for $|t| \leq 1$, we have $h (t) = ( r - t)^2 \li ( - \f{r
t}{1 - r} - \f{t^2}{2} + \f{t^3}{3} - \f{t^4}{4} + \cd \ri )$ and \bee  h(t) + h(-t) & = & \li ( \f{4 r^2}{1 - r} - r^2 \ri ) t^2 -
2 \sum_{k = 2}^\iy \li [ \f{r^2}{2 k} + \f{2r }{2 k - 1} + \f{1}{2 (k-1)} \ri ] t^{2 k}\\
&  > & 2 t^2 \li [ \f{2 r^2}{1 - r} - \f{r^2}{2}  - \li ( \f{r^2}{4} + \f{2 r}{3} + \f{1}{2} \ri ) \f{t^2}{1 - t^2} \ri ]\\
&  \geq & 2 t^2 \li [ \f{2 r^2}{1 - r} - \f{r^2}{2}  - \li ( \f{r^2}{4} + \f{2 r}{3} + \f{1}{2} \ri ) \f{(r \vep)^2}{1 - (r \vep)^2} \ri ]\\
& = & \f{ r^2 t^2  [ 18 + 6 r - (2 r + 6) \vep^2 - 13 (r \vep)^2 - 3 r^3 \vep^2   ] }{ 6 (1 - r) [1 - (r \vep)^2] }\\
& \geq & \f{ r^2 t^2  [ 18 + 6 r - (2 r + 6) - 13 r - 3 r   ] }{ 6 (1 - r) [1 - (r \vep)^2] } = \f{ 2 r^2 t^2  } { 1 - (r \vep)^2  } \geq 0.
\eee This shows $\f{ \pa  \mscr{M}_{\mrm{P}} \li ( z, \se \ri )  }{\pa \vep} - \f{ \pa \mscr{M}_{\mrm{P}} \li ( w, \vse \ri )  }{\pa \vep} > 0$.
The proof of the lemma can be completed by making use of this result and the observation that $\mscr{M}_{\mrm{P}} ( z, \se ) =
\mscr{M}_{\mrm{P}} ( w, \vse)$ for $\vep = 0$.

\epf

\beL

\la{vip89posrev}

{\small \[ \li \{ \f{\ovl{X}_n}{1 + \vep}  \leq  L_n^\ell \leq U_n^\ell \leq \f{\ovl{X}_n}{1 - \vep} \ri \} \supseteq  \li \{ \ovl{X}_n > 0, \;
\mscr{M}_{\mrm{P}} \li ( \f{\ovl{X}_n}{1 + \vep} \li ( 1 + \f{n \vep}{n \vee m_\ell} \ri ), \f{\ovl{X}_n}{1 + \vep} \ri ) \leq \f{1}{m_\ell} \ln
\f{\de_\ell}{2} \ri \}
\]} for all $n \in \mscr{N}$ and $\ell \in \bb{N}$.

\eeL

\bpf  From the definitions of $\mscr{M}_{\mrm{P}}$ and $L_n^\ell$, it is clear that \be \la{comaposrev} \li \{ \mscr{M}_{\mrm{P}} \li ( \f{
\ovl{X}_n}{1 + \vep} \li ( 1 + \f{n \vep}{n \vee m_\ell} \ri ), \f{\ovl{X}_n}{1 + \vep} \ri ) \leq \f{1}{m_\ell} \ln \f{\de_\ell}{2}, \;
\ovl{X}_n = 0 \ri \} = \li \{ \ovl{X}_n = 0, \; \f{\ovl{X}_n}{ 1 + \vep } \leq L_n^\ell \ri \} \ee for $n \in \mscr{N}$.  By Lemma \ref{lempos}
and the definition of $L_n^\ell$, \be \la{combposrev} \li \{ \mscr{M}_{\mrm{P}} \li ( \f{ \ovl{X}_n}{1 + \vep} \li ( 1 + \f{n \vep}{n \vee
m_\ell} \ri ), \f{\ovl{X}_n}{1 + \vep} \ri ) \leq \f{1}{m_\ell} \ln \f{\de_\ell}{2}, \; \ovl{X}_n
> 0 \ri \} = \{ 0 < \ovl{X}_n  \leq L_n^\ell \} \ee for $n \in \mscr{N}$.  It follows from  (\ref{comaposrev}) and (\ref{combposrev}) that  {\small
\be \la{use899a}  \li \{ \f{\ovl{X}_n}{1 + \vep}  \leq  L_n^\ell  \ri \} = \li \{ \mscr{M}_{\mrm{P}} \li ( \f{\ovl{X}_n}{1 + \vep} \li ( 1 +
\f{n \vep}{n \vee m_\ell} \ri ), \f{\ovl{X}_n}{1 + \vep} \ri ) \leq \f{1}{m_\ell} \ln \f{\de_\ell}{2} \ri \} \ee} for all $n \in \mscr{N}$. From
the definitions of $\mscr{M}_{\mrm{P}}$ and $U_n^\ell$, it is clear that \be \la{use899b}  \li \{ \ovl{X}_n > 0, \; \mscr{M}_{\mrm{P}} \li (
\f{\ovl{X}_n}{1 - \vep} \li ( 1 - \f{n \vep}{n \vee m_\ell} \ri ), \f{\ovl{X}_n}{1 - \vep} \ri ) \leq \f{1}{m_\ell} \ln \f{\de_\ell}{2} \ri \}
\subseteq \li \{ \f{\ovl{X}_n}{ 1 - \vep } \geq U_n^\ell \ri \} \ee for $n \in \mscr{N}$.  Finally, combing (\ref{use899a}), (\ref{use899b}) and
using Lemma \ref{vip898posrev} complete the proof of the lemma.

\epf

\beL \la{FinitestoprevPos} $\Pr \{ \f{\ovl{X}_n}{1 + \vep} \leq \bs{\mcal{L}}_n \leq \bs{\mcal{U}}_n \leq \f{\ovl{X}_n}{1 - \vep} \; \tx{for
some} \;  n \in \mscr{N} \} = 1$.

\eeL

\bpf

By the definition of $\bs{\mcal{L}}_n$ and $\bs{\mcal{U}}_n$, it
suffices to show that {\small $\Pr \{ \f{\ovl{X}_{m_\ell} }{1 +
\vep} \leq L_{m_\ell}^\ell \leq U_{m_\ell}^\ell \leq
\f{\ovl{X}_{m_\ell}}{1 - \vep} \; \tx{for some} \;  \ell \in \bb{N}
\} = 1$}.  From Lemma \ref{vip89posrev}, it can be seen that \bee &
& \Pr \li \{ \f{\ovl{X}_{m_\ell} }{1 + \vep} \leq L_{m_\ell}^\ell
\leq U_{m_\ell}^\ell \leq \f{\ovl{X}_{m_\ell}}{1 - \vep} \; \tx{for
some}
\;  \ell \in \bb{N} \ri \}\\
&   & \geq \Pr \li \{ \ovl{X}_{m_\ell} > 0, \; \mscr{M}_{\mrm{P}}
\li ( \ovl{X}_{m_\ell}, \f{\ovl{X}_{m_\ell}}{1 + \vep} \ri ) \leq
\f{1}{m_\ell} \ln \f{\de_\ell}{2}  \; \tx{for some} \; \ell \in
\bb{N} \ri \}. \eee  This inequality and Bonferroni's inequality
imply that \bee &  &  \Pr \li \{ \f{\ovl{X}_{m_\ell} }{1 + \vep}
\leq L_{m_\ell}^\ell \leq U_{m_\ell}^\ell \leq
\f{\ovl{X}_{m_\ell}}{1 - \vep} \; \tx{for some}
\;  \ell \in \bb{N} \ri \}\\
&   & \geq \lim_{\ell \to \iy} \Pr \{ \ovl{X}_{m_\ell}
> 0 \} + \lim_{\ell \to \iy} \Pr \li \{ \mscr{M}_{\mrm{P}}
\li ( \ovl{X}_{m_\ell}, \f{\ovl{X}_{m_\ell}}{1 + \vep} \ri ) \leq \f{1}{m_\ell} \ln \f{\de_\ell}{2} \ri \} - 1. \eee Since $\lm > 0$, it follows
from the law of large numbers that $\lim_{\ell \to \iy} \Pr \{ \ovl{X}_{m_\ell}
> 0 \} = 1$.   To complete the proof of the lemma, it remains  to show that
$\lim_{\ell \to \iy} \Pr \li \{ \mscr{M}_{\mrm{P}} \li ( \ovl{X}_{m_\ell}, \f{\ovl{X}_{m_\ell}}{1 + \vep} \ri ) \leq \f{1}{m_\ell} \ln
\f{\de_\ell}{2} \ri \} = 1$.  This is accomplished as follows.

Let $0 < \eta < 1$.  Noting that $\f{1}{m_\ell} \ln \f{\de_\ell}{2}$
is negative for any $\ell
> 0$ and that {\small $\f{1}{m_\ell} \ln \f{\de_\ell}{2} \to 0
> \mscr{M}_{\mrm{P}} ( \eta \lm, \f{\eta \lm}{1 + \vep} )$} as $\ell
\to \iy$, we have that there exists an integer $\ka$ such that
{\small $\mscr{M}_{\mrm{P}} ( \eta \lm, \f{\eta \lm}{1 + \vep} ) <
\f{1}{m_\ell} \ln \f{\de_\ell}{2}$} for all $\ell \geq \ka$. For
$\ell$ no less than such $\ka$, we claim that $z < \eta \lm$ if
{\small $\mscr{M}_{\mrm{P}} ( z, \f{z}{1 + \vep} ) > \f{1}{m_\ell}
\ln \f{\de_\ell}{2}$} and $z \in [0, \iy)$. To prove this claim,
suppose, to get a contradiction, that $z \geq \eta \lm$. Then, since
{\small $\mscr{M}_{\mrm{P}} ( z, \f{z}{1 + \vep} )$} is
monotonically decreasing with respect to $z \in (0, \iy)$, we have
{\small $\mscr{M}_{\mrm{P}} (z, \f{z}{1 + \vep}  ) \leq
\mscr{M}_{\mrm{P}} ( \eta \lm, \f{ \eta \lm}{1 + \vep} ) <
\f{1}{m_\ell} \ln \f{\de_\ell}{2}$}, which is a contradiction.
Therefore, we have shown the claim and it follows that {\small $ \{
\mscr{M}_{\mrm{P}} ( \ovl{X}_{m_\ell}, \f{\ovl{X}_{m_\ell}}{ 1 +
\vep } )
> \f{1}{m_\ell} \ln \f{\de_\ell}{2} \} \subseteq \{ \ovl{X}_{m_\ell} < \eta
\lm  \}$} for $\ell \geq \ka$. So, {\small
\[  \Pr \li \{ \mscr{M}_{\mrm{P}}
\li ( \ovl{X}_{m_\ell}, \f{\ovl{X}_{m_\ell}}{ 1 + \vep } \ri )
> \f{1}{m_\ell} \ln \f{\de_\ell}{2} \ri \} \leq  \Pr \{ \ovl{X}_{m_\ell}  <
\eta \lm  \}
\]} for large enough $\ell$.  By the law of large numbers,
$\Pr \{ \ovl{X}_{m_\ell} < \eta \lm  \} \to 0$  for sufficiently large $\ell$.  Thus, $\lim_{\ell \to \iy} \Pr \li \{ \mscr{M}_{\mrm{P}} \li (
\ovl{X}_{m_\ell}, \f{\ovl{X}_{m_\ell}}{1 + \vep} \ri ) \leq \f{1}{m_\ell} \ln \f{\de_\ell}{2} \ri \} = 1$.  This proves the lemma.

\epf

\bsk

Now we are in a position to prove that stopping rule F ensures the desired level of coverage probability.  From Lemma \ref{vip89posrev}, we know
that the stopping rule implies that ``continue sampling until $\{ \f{\ovl{X}_n}{1 + \vep}  \leq  L_n^\ell \leq U_n^\ell \leq \f{\ovl{X}_n}{1 -
\vep} \}$ for some $\ell \in \bb{N}$ and $n \in \mscr{N}$''.  We claim that  this stopping rule implies that ``continue sampling until $\{
\f{\ovl{X}_n}{1 + \vep} \leq \bs{\mcal{L}}_n \leq \bs{\mcal{U}}_n \leq \f{\ovl{X}_n}{1 - \vep} \}$ for some $n \in \mscr{N}$''. To show this
claim, we need to show
\[
\bigcup_{\ell \in \bb{N}} \bigcup_{n \in \mscr{N}} \li \{ \f{\ovl{X}_n}{ 1 + \vep }  \leq  L_n^\ell \leq U_n^\ell \leq \f{\ovl{X}_n}{1 - \vep}
\ri \} \subseteq \bigcup_{n \in \mscr{N}} \li \{ \f{\ovl{X}_n}{1 + \vep} \leq \bs{\mcal{L}}_n \leq \bs{\mcal{U}}_n \leq \f{\ovl{X}_n}{1 - \vep}
\ri \},
\]
which follows from the fact that {\small $\bigcup_{\ell \in \bb{N}}  \li \{ \f{\ovl{X}_n} {1 + \vep}  \leq  L_n^\ell \leq U_n^\ell \leq
\f{\ovl{X}_n}{1 - \vep} \ri \} \subseteq \li \{ \f{\ovl{X}_n}{ 1 + \vep} \leq \bs{\mcal{L}}_n \leq \bs{\mcal{U}}_n \leq \f{\ovl{X}_n}{1 - \vep}
\ri \}$}  for every $n \in \mscr{N}$.  From Lemma \ref{FinitestoprevPos}, we know that the sampling process will eventually terminate. It
follows from Lemma \ref{BonferoonrevPos} and Theorem \ref{General Inclusion Principle} that $\Pr \{ | \ovl{X}_{\mbf{n}} - \lm | < \vep \lm \}
\geq 1 - \de$.

\sect{Proof of Theorem \ref{seegood}} \la{seegood_app}

To show (\ref{good1}), note that {\small \bee [( 1 - \mu )^2 + \se]^2 \f{\pa \psi (z, \mu, \se)}{\pa \mu}  = 2 (\mu - z) [( 1 - \mu )^2 + \se] +
(1 - z) [\se - (1 - \mu)^2 ] \ln \li ( 1 + \f{(z - \mu) [ \se + ( 1 - \mu )^2 ] } {\se (1 - z) } \ri ). \eee} If $\se  < (1 - \mu)^2$, then
$\f{\pa \psi (z, \mu, \se)}{\pa \mu} < 0$. If $\se  > (1 - \mu)^2$, then {\small \bee \f{ [( 1 - \mu )^2 + \se]^2 }{ (1 - z) [\se  - (1 - \mu)^2
] } \f{\pa \psi (z, \mu, \se)}{\pa \mu} & = & \ln \li ( 1 + \f{(z - \mu) [ \se + ( 1 - \mu )^2 ] } {\se (1 - z) } \ri ) - \f{ 2 (z - \mu) [( 1 -
\mu )^2 + \se] }{(1 - z) [\se  - (1 -
\mu)^2 ]}  \\
&  \leq & \f{(z - \mu) [ \se + ( 1 - \mu )^2 ] } {\se (1 - z) } - \f{ 2 (z - \mu) [( 1 - \mu )^2 + \se] }{(1 - z) [\se  - (1 - \mu)^2
]}  \\
&  \leq &  - \f{(z - \mu) [ \se + ( 1 - \mu )^2 ] } {\se (1 - z) } \leq 0  \eee} for $0 < \mu < z$.

To show (\ref{good2}),  note that \bee &  & \f{\pa \psi (z, \mu, \se)}{\pa \se}  =  \f{(1 - \mu) (1 - z)}{[ (1 - \mu)^2 + \se ]^2} \li [ \ln \li
( 1 + \f{ (z - \mu) [ (1 - \mu)^2 + \se ] }{\se (1 - z)}  \ri ) - \f{ (z - \mu) [ (1 - \mu)^2 + \se) }{\se (1 - z)} \ri ] \leq 0. \eee

To show (\ref{good3}),  note that {\small \bee ( \mu^2 + \se )^2 \f{\pa \varphi (z, \mu, \se)}{\pa \mu}  =  2 (\mu - z) ( \mu^2 + \se )  + z (
\mu^2 - \se ) \ln \li ( 1 + \f{(\mu - z) ( \se + \mu^2 ) } {\se z } \ri ). \eee} If $\se < \mu^2$, then $\f{\pa \varphi (z, \mu, \se)}{\pa \mu}
> 0$. If $\se
> \mu^2$, then {\small \bee \f{ ( \mu^2 + \se )^2 }{ z (\se - \mu^2 ) } \f{\pa \varphi (z, \mu, \se)}{\pa \mu} & =  & \f{ 2 (\mu - z) ( \mu^2 + \se
) }{z (\se  - \mu^2 )} - \ln \li ( 1 + \f{(\mu - z) ( \se + \mu^2 ) }
{\se z } \ri )  \\
&  \geq & \f{ 2 (\mu - z) ( \mu^2 + \se ) }{z (\se  - \mu^2)} - \f{(\mu - z) [ \se + \mu^2 ] } {\se z } \geq \f{(\mu - z) ( \se + \mu^2 ) } {\se
z } \geq 0  \eee} for $0 < z < \mu$.

To show (\ref{good4}), note that \bee &  &  \f{ \pa \varphi (z, \mu, \se)}{\pa \se}
 = \f{z \mu} {(\mu^2 + \se)^2} \li [ \ln \li ( 1 + \f{ (\mu - z) (\mu^2 + \se) }
 { z \se} \ri )  - \f{ (\mu - z) (\mu^2 + \se) } { z \se} \ri ]  \leq 0 \eee for $0 < z < \mu$.  This completes the proof of the theorem.

 \sect{Proof of Theorem \ref{meaninterval}} \la{meaninterval_app}

Define $L ( \ovl{X}, \se ) = \inf  \li \{ \nu \in [0, \ovl{X}]: \psi ( \ovl{X}, \nu, \se) < \f{\ln \f{3}{\de} }{n} \ri \}$ and $U ( \ovl{X}, \se
) = \sup \li \{ \nu \in [ \ovl{X}, 1 ]: \varphi ( \ovl{X}, \nu, \se) < \f{\ln \f{3}{\de} }{n} \ri \}$ for $0 < \se \leq \f{1}{4}$.  For
simplicity of notations, let $G_{\ovl{X}} (z) = \Pr \{ \ovl{X} \geq z \}$. From the definition of $L ( \ovl{X}, \se )$, Hoeffding's inequality,
and (\ref{good1}) of Theorem \ref{seegood}, we have \bee &  &  \{ L ( \ovl{X}, \se) \geq \mu \} \subseteq \li \{ \mu \leq
\ovl{X}, \;  \psi ( \ovl{X}, \mu, \se) \geq \f{\ln \f{3}{\de} }{n} \ri \}\\
&  &  = \li \{ \mu \leq \ovl{X}, \;  \psi ( \ovl{X}, \mu, \se) \geq \f{\ln \f{3}{\de} }{n}, \;  G_{\ovl{X}} (\ovl{X}) \leq \exp ( - n \psi (
\ovl{X}, \mu, \se) ) \ri \} \subseteq \li \{   G_{\ovl{X}} (\ovl{X}) \leq \f{\de}{3} \ri \}.  \eee  It follows that $\Pr  \{ L ( \ovl{X}, \se)
\geq \mu \} \leq \Pr \li \{   G_{\ovl{X}} (\ovl{X}) \leq \f{\de}{3} \ri \} \leq \f{\de}{3}$.  In a similar manner, we can show that $\Pr  \{ U (
\ovl{X}, \se) \leq \mu \} \leq  \f{\de}{3}$.  Therefore, by Bonferroni's inequality, \be \la{nowaa}
 \Pr  \{ L ( \ovl{X}, \se) < \mu <  U ( \ovl{X}, \se)  \} \geq 1 - \f{2 \de}{3}. \ee Now define {\small $\mcal{U} ( \ovl{X}, \ovl{V}, \mu ) = \sup \li \{  \vse \in \li [ W_\mu,
\f{1}{4} \ri ]: \phi ( W_\mu, \vse ) < \f{\ln \f{3}{\de} }{n} \ri \}$} for $0 < \mu < 1$. Let $ Y = (X - \mu)^2 $.  Then, $Y$ is a random
variable with mean $\se = \si^2$.   Define $\ovl{Y} = \f{\sum_{i=1}^n (X_i - \mu)^2}{n}$.  It can be checked that $\ovl{Y} = W_\mu$.  For
simplicity of notations, let $F_{\ovl{Y}} (y) = \Pr \{ \ovl{Y} \leq y \}$. By Hoeffding's inequality, we have $\Pr \{ \ovl{Y} \leq y \} \leq
\exp \li ( - n \phi (y, \se) \ri )$.  From the definition of $\mcal{U} ( \ovl{X}, \ovl{V}, \mu )$, Hoeffding's inequality,  and the fact that
$\phi ( W_\mu, \se)$ is non-decreasing with respect to $\se > W_\mu$, we have \bee & & \{ \mcal{U} ( \ovl{X},
\ovl{V}, \mu ) \leq \se \} \subseteq \li \{ W_\mu \leq \se, \;  \phi ( W_\mu, \se) \geq \f{\ln \f{3}{\de} }{n} \ri \}\\
&  &  = \li \{ W_\mu \leq \se, \;  \phi ( W_\mu, \se) \geq \f{\ln \f{3}{\de} }{n}, \; F_{\ovl{Y}} (\ovl{Y}) \leq \exp ( - n \phi (W_\mu, \se) )
\ri \} \subseteq \li \{   F_{\ovl{Y}} (\ovl{Y}) \leq \f{\de}{3} \ri \}.  \eee It follows that \be \la{nowbb}
 \Pr  \{ \mcal{U} ( \ovl{X},
\ovl{V}, \mu  ) \leq \se \} \leq \Pr \li \{   F_{\ovl{Y}} (\ovl{Y}) \leq \f{\de}{3} \ri \} \leq \f{\de}{3}. \ee Define random region {\small \[
\mcal{D} ( \ovl{X}, \ovl{V}) = \li \{ (\nu, \vse): 0 < \nu < 1, \; 0 < \vse \leq \nu (1 - \nu), \; L ( \ovl{X}, \vse ) < \nu < U ( \ovl{X}, \vse
), \; \mcal{U} ( \ovl{X}, \ovl{V}, \nu ) > \vse \ri \}. \]} Then, from Bonferroni's inequality, $\Pr \{ (\mu, \se) \in \mcal{D} ( \ovl{X},
\ovl{V} ) \} \geq 1 - \de$. By (\ref{good1}), (\ref{good3}) of Theorem \ref{seegood}, and the fact that $\phi ( W_\nu, \vse)$ is non-decreasing
with respect to $\vse > W_\nu$, we have  $\mcal{D} ( \ovl{X}, \ovl{V}) = \mcal{A} \cup \mcal{B}$, where \bee &  & \mcal{A} = \li \{ (\nu, \vse):
\nu \in (0, \ovl{X}), \; \vse \in \li (0, \nu (1 - \nu) \ri ], \; \max \li [  \psi ( \ovl{X}, \nu,
\vse), \; \phi (W_\nu, \vse) \; \bb{I}_{\{\vse > W_\nu\}}  \ri ] <  \f{\ln \f{3}{\de}}{n}   \ri \},\\
&  & \mcal{B} =  \li \{ (\nu, \vse): \nu \in (\ovl{X}, 1),  \; \vse \in \li (0, \nu (1 - \nu) \ri ], \; \max \li [  \varphi ( \ovl{X}, \nu,
\vse), \; \phi (W_\nu, \vse ) \; \bb{I}_{\{\vse > W_\nu\}}  \ri ] < \f{\ln \f{3}{\de}}{n}   \ri \}. \eee Finally, the theorem follows from the
observation that $\mcal{D} ( \ovl{X}, \ovl{V}) = \mcal{A} \cup \mcal{B}$.

\sect{Proof of Theorem \ref{Region}} \la{Region_app}

Define {\small \bee  L ( \ovl{X}, \se ) = \inf  \li \{ \nu \in [0, \ovl{X}]: \psi ( \ovl{X}, \nu, \se) < \f{\ln \f{4}{\de} }{n} \ri \}, \qqu U (
\ovl{X}, \se ) = \sup \li \{ \nu \in [ \ovl{X}, 1 ]: \varphi ( \ovl{X}, \nu, \se) < \f{\ln \f{4}{\de} }{n} \ri \} \eee} for $0 < \se \leq
\f{1}{4}$.  Define {\small \bee  \mcal{L} ( \ovl{X}, \ovl{V}, \mu ) = \inf \li \{  \vse \in [0, W_\mu]: \phi ( W_\mu, \vse) < \f{\ln \f{4}{\de}
}{n} \ri \}, \qqu  \mcal{U} ( \ovl{X}, \ovl{V}, \mu ) = \sup \li \{  \vse \in \li [ W_\mu, \f{1}{4} \ri ]: \phi ( W_\mu, \vse ) < \f{\ln
\f{4}{\de} }{n} \ri \} \eee} for $0 < \mu < 1$.   By a similar method as that for proving (\ref{nowaa}), we can show that $\Pr  \{ L ( \ovl{X},
\se) < \mu < U ( \ovl{X}, \se)  \} \geq 1 - \f{\de}{2}$.  By a similar method as that for proving (\ref{nowbb}), we can show that
\[
\Pr  \{ \mcal{L} ( \ovl{X}, \ovl{V}, \mu  ) \geq \se \} \leq \f{\de}{4}, \qqu \Pr  \{ \mcal{U} ( \ovl{X}, \ovl{V}, \mu  ) \leq \se \} \leq
\f{\de}{4}.
\]
Therefore, by Bonferroni's inequality, $\Pr  \{ \mcal{L} ( \ovl{X}, \ovl{V}, \mu  ) < \se <  \mcal{U} ( \ovl{X}, \ovl{V}, \mu  ) \} \geq 1 -
\f{\de}{2}$.  Again by Bonferroni's inequality,
\[
\Pr \{  L ( \ovl{X}, \se ) < \mu < U ( \ovl{X}, \se ), \; \mcal{L} ( \ovl{X}, \ovl{V}, \mu ) < \se < \mcal{U} ( \ovl{X}, \ovl{V}, \mu ) \} \geq
1 - \de.
\]
By (\ref{good1}), (\ref{good3}) of Theorem \ref{seegood} and the unimodal property of $- \phi (W_\mu, \se)$ with respect to $\se$, we have  that
{\small
\[ \mscr{D} ( \ovl{X}, \ovl{V}) = \li \{ (\nu, \vse): 0 < \nu < 1, \; 0 < \vse \leq \nu (1 - \nu), \; L ( \ovl{X}, \vse ) < \nu < U ( \ovl{X},
\vse ), \; \mcal{L} ( \ovl{X}, \ovl{V}, \nu ) < \vse < \mcal{U} ( \ovl{X}, \ovl{V}, \nu ) \ri \},
\]}
which implies that $\Pr \{ (\mu, \si^2) \in \mscr{D} ( \ovl{X}, \ovl{V} ) \} \geq 1 - \de$.  This completes the proof of the theorem.

\end{document}